\documentclass[preprint,3p]{elsearticle}
\usepackage{graphicx} 
\usepackage{algorithm}
\usepackage{amsmath}
\usepackage{amsfonts} 
\usepackage{xcolor} 
\usepackage{algpseudocode}
\usepackage{tabularx} 
\usepackage{bbm}
\usepackage{amssymb}
\newcommand{\PP}[0]{\mathbb{P}}
\newcommand{\setR}[0]{\mathbb{R}}
\newcommand{\Fhat}[0]{\widehat{F}}

\newcommand{\E}[0]{\mathbb{E}}
\newcommand{\LQ}[0]{\left[}
\newcommand{\D}[0]{\mathcal{D}}
\newcommand{\RQ}[0]{\right]}

\newcommand{\xb}[0]{\mathbf{x}}
\newcommand{\xib}[0]{\boldsymbol\xi}
\newcommand{\argmin}[0]{\text{argmin}}
\newcommand{\Cvar}[0]{\text{CVaR}_\beta}
\newcommand{\Var}[0]{\text{VaR}_\beta}
\newcommand{\Variance}[0]{\text{Var}}

\newcommand{\Tol}[0]{\text{Tol}}
\newcommand{\I}[0]{\mathcal{I}}
\newcommand{\RB}[0]{\text{RB}}
\newcommand{\C}[0]{\mathcal{C}}
\newcommand{\rhow}[0]{\widetilde{\rho}}
\newcommand{\supp}[0]{\text{supp}}
\newcommand{\that}[0]{\widehat{t}}
\newcommand{\ISt}[0]{\text{IS}}
\newcommand{\Hc}[0]{\mathcal{H}}
\newcommand{\dprime}[0]{{\prime \prime}}
\newcommand{\Ftilde}[0]{\widetilde{F}}

\newtheorem{thm}{Theorem}
\newtheorem{prop}[thm]{Proposition}
\newtheorem{lem}[thm]{Lemma}
\newdefinition{rmk}{Remark}
\newproof{pf}{Proof}
\newtheorem{ass}{Assumption}

\makeatletter
\def\BState{\State\hskip-\ALG@thistlm}
\makeatother

\begin{document}
\begin{frontmatter}
\title{An adaptive importance sampling algorithm for risk-averse optimization}
\date{\today}
\author[1]{Sandra Pieraccini}
\ead{sandra.pieraccini@polito.it}
\author[1]{Tommaso Vanzan\corref{cor1}}
\ead{tommaso.vanzan@polito.it}
\affiliation[1]{organization={Dipartimento di Scienze Matematiche, Politecnico di Torino},
addressline={Corso Duca degli Abruzzi 24},
postcode={10129},
city={Torino},
country={Italy}}
\cortext[cor1]{Corresponding author}
\begin{abstract}
Adaptive sampling algorithms are modern and efficient methods that dynamically adjust the sample size throughout the optimization process. However, they may encounter difficulties in risk-averse settings, particularly due to the challenge of accurately sampling from the tails of the underlying distribution of random inputs. This often leads to a much faster growth of the sample size compared to risk-neutral problems.
In this work, we propose a novel adaptive sampling algorithm that adapts both the sample size and the sampling distribution at each iteration. The biasing distributions are constructed on the fly, leveraging a reduced-order model of the objective function to be minimized, and are designed to oversample a so-called risk region. As a result, a reduction of the variance of the gradients is achieved, which permits to use fewer samples per iteration compared to a standard algorithm, while still preserving the asymptotic convergence rate.
Our focus is on the minimization of the Conditional Value-at-Risk (CVaR), and we establish the convergence of the proposed computational framework. Numerical experiments confirm the substantial computational savings achieved by our approach.
\end{abstract}
\end{frontmatter}

\section{Introduction}
In this manuscript, we consider the minimization of
\begin{equation}\label{eq:model}
    \min_{z\in \C} \left\{\Fhat(z):= \Cvar \LQ f(z,\xib)\RQ\right\},
\end{equation}
where $z$ is a decision variable, $\C$ is a convex and closed subset of an Hilbert space $(\mathcal{Z},(\cdot,\cdot))$, $\xib$ is a random variable over a probability space $(\Omega,\mathcal{F},\PP)$, with values on a subset $\Gamma\subset \mathbb{R}^N$, $N\in\mathbb{N}^+,$ and with probability density $\rho$. The quantity of interest is $f:\mathcal{Z}\times \Gamma\rightarrow \setR$ which is supposed to be smooth and whose gradient with respect to $z$ is denoted by $\nabla f(z,\xib)$. 
$\Cvar$ is the so-called Conditional Value-at-Risk (CVaR), namely the expectation of $f(z,\cdot)$ conditioned on being above its $\beta$ quantile.
If $\Cvar$ is replaced by the expectation operator $\E_\rho$, \eqref{eq:model} has found numerous applications in different fields such as, e.g., machine learning \cite{bottou2018optimization,Goodfellow-et-al-2016}, portfolio optimization \cite{brandimarte2017introduction,fabozzi2007robust}, stochastic programming \cite{birge2011introduction,shapiro2021lectures} and PDE-constrained optimization with random inputs \cite{antil2018frontiers,martinez2018optimal,Nobile24CT}. More recently, minimization of the $\Cvar$ has drawn interest since one often prefers to be robust against risk-averse events; see, e.g., \cite{rockafellar2000optimization,rockafellar2013fundamental,kouri2018existence,curi2020adaptive}.

The numerical solution of \eqref{eq:model} with $\E_\rho$ has traditionally been dealt with two main classes of algorithms. The first family, called Stochastic Approximation (SA) methods \cite[Chapter 5.9]{shapiro2021lectures}, includes iterative methods that at each iteration draw new realizations of $\xib$ independent from the previous ones and perform one step of a descent algorithm. Examples are the stochastic gradient method and its variants; see, e.g., \cite{bottou2018optimization} and references therein. These algorithms show a low cost per iteration and, under suitable assumption, converge to the minimizer $z^\star$ of \eqref{eq:model}, but at a sublinear rate.
The second family, called Sample Average Approximation (SAA) in stochastic optimization communities or simply batch methods in machine learning, consists in replacing $\PP$, once and for all, with a discrete empirical probability distribution $\PP_M$, obtained by sampling $M$ realizations of $\xib$. The approximated problem can then be solved with standard nonlinear optimization algorithms that often converge very rapidly, but only to a minimizer $z^\star_M\neq z^\star$.
Between these two classes, there are attractive modern strategies such as adaptive sampling methods (see, e.g., \cite{byrd2012sample,friedlander2012hybrid,bollapragada2018adaptive,beiser2023adaptive,xie2024constrained}) which aim to improve the efficiency of SA methods by dynamically adjusting the sample size along the optimization. The driving idea is to use cheap, but
noisy, gradient evaluations far from the optimum, while more accurate and expensive estimates are computed close to the minimizer. Earlier works showed that increasing the sample size geometrically is sufficient to preserve linear convergence of a gradient descent algorithm, while a few implementable rules have been now derived to \textit{adapt} the sample size based on current optimization progresses, both in the unconstrained and constrained settings.

Compared to $\E_\rho$, both the evaluation and optimization of $\Cvar$ raise additional challenges. First, evaluating $\Cvar\LQ f(z,\xib)\RQ$ requires to sample from the tail of the unknown distribution of $f(z,\cdot)$, which implies that a large number of samples are needed to have sufficiently accurate estimates. Second, $\nabla f(z,\xib)$ actually contributes to the approximation of $\nabla \Fhat(z)$ only for samples $\xib$ belonging to a so-called \textit{risk-region}, defined as the set of $\xib$ such that $f(z,\xib)\geq \Var\LQ f(z,\xib)\RQ$, $\Var\LQ f(z,\xib)\RQ$ being the Value-at-Risk of level $\beta$, that is, its $\beta$ quantile. Thus a significant computational effort is wasted if evaluating $f(z,\xib)$ is expensive, as one may compute $f(z,\xib)$ and then realize that the computation of $\nabla f(z,\xib)$ is not needed. This is typically the case in PDE-constrained optimization under uncertainty, where evaluating $f(z,\xib)$ consists in solving a PDE state equation, while $\nabla f(z,\xib)$ entails the solution of an additional PDE adjoint equation.
Third, adaptive sampling algorithms commonly adjust the sample size along the optimization to control the variance of $\xib\rightarrow \nabla f(z,\xib)$. However, this variance typically grows as $\beta$ increases, leading to a faster increase of the sample sizes, hence making risk-averse optimization significantly more expensive compared to the risk-neutral case.

\subsection{Contributions and related works}
The primary contribution of this manuscript is the derivation of a new adaptive sampling algorithm that adjusts the sampling distribution along the optimization. The algorithm generates a
sequence of biasing densities $\left\{\widetilde{\rho}_k\right\}_{k\in \mathbb{N}}$ that are adapted to the sequence of iterates $\left\{z_k\right\}_{k\in \mathbb{N}}$, so that each $\widetilde{\rho}_k$ oversamples the current risk-region, and thus most of the sample do contain valuable gradient information.
This leads to an importance sampling estimation of the gradient $\nabla \Fhat(z)$, and to a variance reduction which permits then to take smaller sample sizes compared to standard adaptive sampling algorithms, eventually reducing the overall computational cost.
Our algorithm \textit{learns} the biasing distributions on the fly during the optimization, relying on the availability of a reduced order model (ROM) of the map $\xib\rightarrow f(z_k,\xib)$, denoted by $f_r(z_k,\cdot)$, at every iterate $z_k$. The ROM is then used both to estimate the current $\Var\LQ f(z_k,\cdot)\RQ$ and to generate the new biased samples using an acceptance-rejection algorithm.
The construction and availability of such a ROM is of course strongly problem dependent. We here present numerical experiments based on an application to PDE-constrained optimization under uncertainty, and discuss how these reduced models can be efficiently built online, during the optimization, at a reduced computational cost.
A secondary contribution of this manuscript lies in a convergence analysis of an alternating optimization algorithm tailored for (1). Indeed, while \cite{beiser2023adaptive} presents an alternating procedure without discussing its convergence, \cite{ganesh2023gradient} presents a rigorous convergence analysis under however some quite restrictive hypothesis. We here propose an alternative path to prove convergence under some relaxed assumptions.

Besides the above mentioned references on adaptive sampling algorithms, we wish to emphasize a couple of manuscripts that are more closely related to this work. In machine learning, \cite{curi2020adaptive} precisely addresses the challenges of the $\Cvar$ optimization using adaptive sampling algorithm, and similarly proposes to adjust the sampling distribution along the optimization. However, the construction of the biasing distribution differs and it is based on a dual reformulation of \eqref{eq:model} as a distributionally robust problem, leading to a min-max game between two players. 
Closely related to this work are also \cite{heinkenschloss2018conditional,heinkenschloss2020adaptive}, which propose two approaches to use ROMs for the evaluation of the $\Cvar$. 
The first approach consists in simply replacing the full order model with a ROM, while the second suggests to use a ROM to derive a biasing distribution from which then sampling the full order model. The first approach has been applied to the $\Cvar$ minimization within SAA algorithms in \cite{zou2022locally,Markowski2022}. In contrast, as the importance sampling approach cannot be embedded into SAA algorithms (since the sample set is fixed in these algorithms), the second approach has been less investigated. The present manuscript, which combines importance sampling with adaptive sampling algorithms, can be viewed as an extension of the second approach of \cite{heinkenschloss2018conditional} to the optimization context.

The manuscript is structured as follows. Section \ref{sec:problem_setting} introduces the model and a classical reformulation of \eqref{eq:model} as a double minimization problem. Then, an alternating optimization procedure is recalled and subsection \ref{sec:challenges} discusses in details the challenges of adaptive sampling algorithms for risk-averse minimization. Section \ref{sec:adaptive_importance_sampling} recalls the importance sampling framework, discusses the ROM-based construction of the biasing densities, and present the adaptive importance sampling algorithm. Subsection \ref{sec:convergence} proposes a convergence analysis of the alternating procedure for the interested readers.
Section \ref{sec:numerical} presents numerical experiments aimed at showing the better complexity of the new approach and details the dynamic update of the ROM during the optimization.

\section{Problem setting}\label{sec:problem_setting}
In this section, we recall an equivalent formulation of \eqref{eq:model} that is the starting point of several approaches to optimize $\Cvar$ \cite{rockafellar2000optimization,kouri2016risk,ciaramella2024multigrid}. We then describe an adaptive sampling algorithm to minimize \eqref{eq:model}, and highlight a few challenges that will be addressed by our contributions in Section \ref{sec:adaptive_importance_sampling}.

\subsection{Conditional Value-At-Risk}
For a given random variable $X\in L_{\rho}^1(\Gamma;\setR)$, the Value-at-Risk of confidence level $0<\beta<1$ (also known as $\beta$ quantile) is defined as
\[\Var\LQ X\RQ:=\inf_{t\in \setR}\left\{\PP(X\leq t)\geq \beta\right\},\]
and its existence is guaranteed by the right continuity of the cumulative distribution function.
The Conditional Value-at-Risk (also called expected shortfall or average Value-at-Risk) of confidence level $\beta$ is the expectation of $X$ conditioned on being larger than $\Var\LQ X\RQ$,
\[\Cvar\LQ X\RQ:=\E_{\rho}\LQ X|X\geq \Var\LQ X\RQ\RQ .\]
The Conditional Value-at-Risk has gained popularity as risk measure since, in contrast to the Value-at-Risk, it embeds information on the average loss above the critical $\beta$ quantile.
In \cite{rockafellar2000optimization,rockafellar2002conditional}, Rockafellar and Uryasev proved that $\Cvar\LQ X\RQ$ admits the equivalent formulation 
\begin{equation}\label{eq:Cvarequidef}
\Cvar\LQ X\RQ= \inf_{t\in \setR}\left\{t+\frac{1}{1-\beta}\E_{\rho}\LQ (X-t)^+\RQ \right\},
\end{equation}
where $(\cdot)^+:=\max(0,\cdot)$, and that the infimum  $\eqref{eq:Cvarequidef}$ is attained, although possibly on a whole bounded interval of the real line $\mathcal{A}_{\beta}\LQ X\RQ$ \footnote{$\mathcal{A}_{\beta}\LQ X\RQ$ is not a singleton when the cumulative distribution function is flat in a right neighborhood of $\Var\LQ X\RQ$, see \cite[Section 2]{rockafellar2002conditional}.}, with $\Var\LQ X\RQ=\min_{t}\left\{t\in \mathcal{A}_{\beta}\LQ X\RQ\right\}$. 

This result allows us to reformulate \eqref{eq:model} as the double minimization
\begin{equation}\label{eq:model_double}
    \min_{z\in \C,\; t\in \setR} \left\{\Ftilde(z,t):=t+\frac{1}{1-\beta}\E_{\rho}\LQ (f(z,\xib)-t)^+\RQ\right\},
\end{equation}
which enjoys attractive features: on the one hand, solving \eqref{eq:model_double} permits to both optimize $\Cvar\LQ f\RQ$ and simultaneously get the associated $\Var\LQ f(z^\star,\cdot)\RQ$, see, e.g., \cite[Corollary 15]{rockafellar2002conditional}. On the other hand, \eqref{eq:model_double} involves the minimization of an expectation, and thus either SAA or SA methods could be readily applied.

Nevertheless \eqref{eq:model_double} involves the function $(\cdot)^+$ which is only Lipschitz continuous. Despite $\Ftilde$ may be jointly Fr\'{e}chet differentiable with respect to $z$ and $t$ due to the smoothing effect of $\E_{\rho}$\footnote{Sufficient conditions are the differentiability of $f$ with respect to $z$ and that the distribution of $f(z,\cdot)$ does not have atoms.} \cite[Theorem 2.1]{ganesh2023gradient},  $\Ftilde$ is non-differentiable whenever the continuous expectation is replaced by an empirical approximation (as in SAA methods), while $(f(z,\xib)-t)^+$ does not admit a classical gradient (which prevents gradient-based SA methods). In order to use optimization algorithms for smooth functions, \cite{kouri2016risk} suggested to replace $(\cdot)^+$ with a smooth approximation $g_{\varepsilon}$. 
In this work, we choose the $C^2$ differentiable function 
\begin{equation*}\label{eq:smoothapproximation}
g_{\varepsilon}=\begin{cases}
0\quad &\text{if } x\leq -\frac{\varepsilon}{2},\\
\frac{(x-\frac{\varepsilon}{2})^3}{\varepsilon^2}-\frac{(x-\frac{\varepsilon}{2})^4}{2\varepsilon^3}\quad &\text{if } x\in(-\frac{\varepsilon}{2},\frac{\varepsilon}{2}),\\
x\quad &\text{if }x\geq \frac{\varepsilon}{2},
\end{cases}
\end{equation*}
and remark two properties that will be used in the following, namely that $\sup_{t\in\setR} |g_\varepsilon^\prime(t)|\leq 1$ and $g_{\varepsilon}^\dprime(t)\geq 0$, for all $t\in \setR$.
From now on, we consider the smoothed version of \eqref{eq:model_double},
\begin{equation}\label{eq:model_double_smoothed}
    \min_{z\in \C,\; t\in \setR} \left\{\Ftilde^{\varepsilon}(x,t):=t+\frac{1}{1-\beta}\E_{\rho}\LQ g_{\varepsilon}(f(z,\xib)-t) \RQ\right\}.
\end{equation}
An useful guide to choose a suitable value of $\varepsilon$ is the error estimate \cite[Theorem 4.13]{kouri2016risk},
\[ (|t^\star - t^\star_{\varepsilon}|^2+\|z^\star - z^\star_{\varepsilon}\|^2)^{\frac{1}2}\leq C \sqrt{\varepsilon},\]
where 
$(z^\star,t^\star)$ and $(z^\star_\varepsilon,t^\star_\varepsilon)$ are  minimizers of \eqref{eq:model_double} and \eqref{eq:model_double_smoothed}, respectively, and $C$ is a positive constant independent of $(z^\star,t^\star)$. Concretely, $\varepsilon$ should be carefully tuned via a tradeoff between a desired relative error, and the conditioning of \eqref{eq:model_double_smoothed} which deteriorates as $\varepsilon\rightarrow 0$.

The double minimization formulation paves the way to a family of algorithms that alternate between exact/inexact minimization with respect to $x$ and $t$, see \cite{Markowski2022,ganesh2023gradient,beiser2023adaptive,ciaramella2024multigrid}. A particular instance is the procedure
\begin{align}
t_{k}&\in \argmin_{t\in \setR} \Ftilde^{\varepsilon}(z_k,t),\label{eq:mini_t}\\
z_{k+1}&=z_k+\alpha_k d_k,\label{eq:mini_x}
\end{align}
where $d_k$ is an appropriate direction, hopefully descent for $\Ftilde^{\varepsilon}$ at $(z_k,t_{k})$.
Notice that \eqref{eq:mini_t} consists in a simple one-dimensional optimization problem which can be solved either using some statistical library to compute quantiles or through, e.g., a bisection algorithm applied to $\nabla_t \Ftilde^{\varepsilon}$. In our numerical implementation, we followed the latter approach. Instead, for a fixed $t$, the direction $d_k$ can be computed with any optimization algorithm mentioned in the introduction, see, e.g., \cite{beiser2023adaptive,ganesh2023gradient} for adaptive sampling approaches which we recall in the next subsection.

\subsection{An adaptive sampling algorithm for risk-neutral optimization}

A traditional mini-batch stochastic gradient algorithm applied to the unconstrained minimization of $F(z):=\E\LQ f(z,\xib)\RQ$ computes at iteration $k=0,\dots$,
\[z_{k+1}=z_k-\alpha_k \frac{1}{M}\sum_{j=1}^M \nabla f(z_k,\xib^k_j),\]
where $\left\{\xib^k_j\right\}_{j=1}^M$, $M\in \mathbb{R}$, are independent realizations sampled according to $\rho$.
This generally leads to noisy gradient estimates, which are then mitigated by a decreasing step size $\alpha_k$, ultimately responsible for the slow convergence.                        
In contrast, adaptive sampling methods reduce the variance of the stochastic gradients by dynamically increase the sample size along the optimization while using a fixed step size.
An iteration reads
\[z_{k+1}=z_k-\alpha \nabla F_{S_k}(z_k)=z_k-\alpha \frac{1}{M_k}\sum_{j=1}^{M_k} \nabla f(z_k,\xib^k_j),\]
where $F_{S_k}(z_k):=\frac{1}{M_k}\sum_{j=1}^{M_k} f(z_k,\xib^k_j)$, $S_k=\left\{\xib^k_j\right\}_{j=1}^{M_k}$ denotes the set of i.i.d. realizations of $\xib$ drawn at iteration $k$, and $\alpha$ is now a fixed-constant suitably tuned. 

Starting from \cite[Theorem 3.1]{friedlander2012hybrid}, which showed that growing the sample size geometrically, i.e. $M_k=\tau^k$ for $\tau>1$, permits to recover the linear convergence of gradient descent with fixed step length for the unconstrained minimization of strongly convex functions, 
several contributions have then developed criteria to adapt the sample size based on the \textit{actual} progress of the unconstrained optimization (and thus not on a preset rule), see \cite{byrd2012sample,hashemi2014adaptive,cartis2018global,bollapragada2018adaptive}. 
In particular, let $\E\LQ\cdot\RQ$ denote the full expectation with respect to the distribution of all sample sets $\left\{S_k\right\}_{k\geq 1}$,
introduce the filtration $\mathcal{F}_k=\sigma(S_i,\;i<k)$, that is the sigma algebra generated by all samples $\left\{\xib^i_j\right\}_{j=1}^{M_i}$ used up to iteration $k-1$, and the shorthand notation $\E_k\LQ \cdot\RQ$ for the conditional expectation $\E\LQ \cdot | \mathcal{F}_k\RQ$. Hence, for instance $\E_k\LQ F_{S_k}(z_k)\RQ$ is well defined and corresponds to the expectation of $F_{S_k}(z_k)$ with respect only to the distribution of the samples $\left\{\xib^k_j\right\}_{j=1}^{M_k}$ drawn at iteration $k$. The authors of
\cite{byrd2012sample} showed that imposing at each iteration that
\begin{equation}\label{eq:norm_test_the}
\E_k\LQ \|\nabla F_{S_k}(z_k)-\nabla F(z_k)\|^2\RQ=  \frac{\E_{\rho}\LQ \|\nabla f(z_k,\cdot)-\nabla F(z_k)\|^2\RQ}{M_k} = \frac{\Variance_{\rho}\LQ \nabla f(z_k,\cdot)\RQ}{M_k} \leq \theta^2 \|\nabla F(z_k)\|^2,
\end{equation}
is sufficient to preserve the linear convergence of the descent algorithm.
In words, \eqref{eq:norm_test_the} requires that the variance of the estimated gradient $\nabla F_{S_k}(z_k)$ decreases proportionally to the norm square of the exact gradient evaluated at the current iteration $z_k$.

In practice, imposing \eqref{eq:norm_test_the} at each iteration may be inconvenient, as it would require to sequentially add new samples to $S_k$ until the criterion is satisfied. Furthermore, the exact gradient is not an available quantity and thus it is approximated with the information available at iteration $k$, namely the empirical gradient. A common approach, often called \textit{norm test}, is then to use \eqref{eq:norm_test_the} to estimate the \textit{next} sample size by computing
\begin{equation}\label{eq:norm_test_pra}
\varrho=\frac{\sum_{\xib^k_j\in S_k}\|\nabla f(z_k,\xib^k_j)-\nabla F_{S_k}(z_k)\|^2}{(M_k-1)(M_k)\theta^2 \|\nabla F_{S_k}(z_k)\|^2} ,
\end{equation}
and set $M_{k+1}=\varrho M_k$ if $\varrho>1$.
The factor $M_k-1$ in \eqref{eq:norm_test_pra} is due to the unbiased estimator of $\Variance_{\rho}\LQ \nabla f(z_k,\cdot)\RQ$.
For convex constrained optimization problems, \eqref{eq:norm_test_the}
is not appropriate since $\nabla F$ may not vanish at the optimum.
Nevertheless, \cite{beiser2023adaptive,xie2024constrained} proved that the alternative condition
\begin{equation}\label{eq:norm_test_pra_con}
\E_k\LQ \|\nabla F_{S_k}(z_k)-\nabla F(z_k)\|^2\RQ\leq \theta^2\|R_{S_k}(z_k)\|^2 
\end{equation}
where $R_{S_k}(z_k):=\frac{1}{\alpha}(z_{k+1}-z_k)$ is sufficient to guarantee linear convergence of the iterates.

\subsection{Challenges associated to the $\Cvar$ optimization}\label{sec:challenges}
In the context of the $\Cvar$ minimization, an adaptive sampling algorithm has been first proposed in \cite{beiser2023adaptive} (see also \cite{ganesh2023gradient} for a slightly different version), and it is recalled in Alg. \ref{alg:adaptive_sampling}. Here, $\mathcal{P}$ denotes the projection operator on the convex set $\C$. Alg. \ref{alg:adaptive_sampling} consists of an alternating optimization procedure. At every iteration $k$, a new set of i.i.d. samples $S_k$ is drawn. Then, the approximated function 
\[\Ftilde_{S_k}^{\varepsilon}(z,t):=t+\frac{1}{(1-\beta)M_k}\sum_{j=1}^{M_k}g_{\varepsilon}(f(z,\xib_j)-t)\]
is first exactly minimized with respect to the $t$ variable, and then a gradient step is performed with respect to $z$. Finally, the norm test (direct extension of \eqref{eq:norm_test_pra})
\begin{equation}\label{eq:norm_test_pra_cvar}
\varrho=\frac{\sum_{\xib^k_j\in S_k}\left\| \frac{g^\prime_{\varepsilon}(f(z_k,\xib_j^k)-t_k)\nabla f(z_k,\xib_j^k)}{1-\beta}-\nabla \Ftilde^{\varepsilon}_{S_k}(z_k,t_k)\right\|^2}{(M_k-1)(M_k)\theta^2 \|R_{S_k}(z_k)\|^2}
\end{equation}
decides whether to increase the sample size. The symbol $\nabla$ still denotes the gradient with respect only to the variable $z$ .
We highlight that at every iteration of Alg. \ref{alg:adaptive_sampling}, $t_{k}$ represents an approximation of $\Var\LQ f(z_k,\cdot)\RQ$ computed with the current $M_k$ samples. This property linking $z_k$ and $t_{k}$ will come at hand in the coming discussion.
\begin{algorithm}
\caption{Adaptive Sampling Algorithm for $\Cvar$ minimization}\label{alg:adaptive_sampling}
\begin{algorithmic}[1]
\State \textbf{Require} $\mathbf{z}^0$, $M_0$, $\alpha$.
\State Set $k=0$.
\State Sample $S_0=\left\{\xib^k_i\right\}_{i=1}^{M_0}$ from $\rho$
\While {Convergence not reached}
\State $t_{k}\in\argmin_{t\in \setR} \Ftilde^{\varepsilon}_{S_k}(z_k,t)$.
\State $z_{k+1}=\mathcal{P}\left( z_k-\alpha_k \nabla \Ftilde^{\varepsilon}_{S_k}(z_k,t_{k})\right)$.
\If {Test \eqref{eq:norm_test_pra_con} is not satisfied}
\State $M_{k+1}=\varrho M_k$, with $\varrho$ as in \eqref{eq:norm_test_pra_cvar}
\State Sample $S_{k+1}=\left\{\xib^{k+1}_i\right\}_{i=1}^{M_{k+1}}$ from $\rho$.
\EndIf
\EndWhile
\end{algorithmic}
\end{algorithm}

In order to illustrate the computational challenges of Alg. \ref{alg:adaptive_sampling}, we start by re-interpreting the quantities that are involved in the $\Cvar$ minimization and in the norm test. We first define the event
\[G_{z,t}:=\left\{\xib:\; f(z,\xib)\geq t\right\},\]
and observe that, if $t=\Var\LQ f(z,\cdot)\RQ$, then $\PP(G_{z,t})=1-\beta$ and the gradient of $\Ftilde^{\varepsilon}$ satisfies
\[\nabla \Ftilde^{\varepsilon}(z,t)=\E\LQ \frac{\nabla g_{\varepsilon}(f(z,\cdot)-t)}{1-\beta}\RQ= \E\LQ \frac{g_{\varepsilon}^\prime(f(z,\cdot)-t) \nabla f(z,\cdot)}{1-\beta}\RQ \approx \E\LQ \nabla f(z,\cdot) | G_{z,t}\RQ,\]
i.e., $\nabla \Ftilde^{\varepsilon}(z,t)$ is an approximation (since $g_{\varepsilon}^\prime(f(z,\xib)-t)$ is a smoothed version of the characteristic function $\mathbbm{1}_{G_{z,t}}$) of the expected value of $\nabla f(z,\cdot)$ conditioned on the event $G_{z,t}$. To approximate $\nabla \Ftilde^{\varepsilon}$ at each iteration, Alg. \ref{alg:adaptive_sampling} considers the function $\Ftilde^{\varepsilon}_{S_k}$ evaluated at the current iterate $(z_k,t_k)$, where $t_k$ is a sample-based approximation of $\Var\LQ f(z_k,\cdot)\RQ$. Consequently, $\nabla \Ftilde^{\varepsilon}(z_k,t_k)$ is estimated by the empirical average
\[\nabla \Ftilde_{S_k}^{\varepsilon}(z_k,t_k)=\frac{1}{(1-\beta)M_k}\sum_{\xib_j^k\in S_k}g_\varepsilon^\prime (f(z_k,\xib_j^k)-t_k)\nabla f(z_k,\xib_j^k).\] 
One first challenge is then due to term $g_\varepsilon^\prime (f(z_k,\cdot)-t_k)$, since only very few samples actually contribute to the gradient estimation, namely those belonging to the so-called risk-region
\begin{equation*}
G^{\varepsilon}_{z_k,t_k}\LQ f\RQ :=\left\{\xib:\; g_{\varepsilon}^\prime(f(z_k,\xib )-t_k)>0 \right \}=\left\{\xib:\; f(z_k,\xib)\geq t_k-\frac{\varepsilon}{2}\right \}.      
\end{equation*}
For small values of $\varepsilon$, approximately a fraction $1-\beta$ of the samples belongs to $G_{z_k,t_k}^{\varepsilon}$ at every iteration, and this leads to a great waste of computational resources to get accurate gradient estimations.

Given this consideration, one might infer that Alg. 1 would require a sequence of sample sizes that are approximately $\frac{1}{1-\beta}$ larger than the sample sizes used to minimize the risk-neutral objective functional $\E\LQ f(z,\cdot)\RQ$. Unfortunately, this is not the case, since Alg.1 sets the sample size according to the \textit{variance} of the gradients, and thus it may lead to much larger values of $M_k$. To understand this, we analyze the theoretical version of the norm test \eqref{eq:norm_test_the},
\begin{equation*}
\begin{aligned}
\Variance_k\LQ \nabla \Ftilde^{\varepsilon}_{S_k}(z_k,t_k)\RQ =\frac{\Variance_{\rho}\LQ \frac{g^\prime_{\varepsilon}(f(z_k,\cdot)-t_k)\nabla f(z_k,\cdot)}{1-\beta}\RQ}{M_k} \leq \theta^2 \|\nabla \Ftilde^{\varepsilon}(z_k,t_k)\|^2,
\end{aligned}
\end{equation*}
and focus on the term $\text{Var}_{\rho}\LQ \frac{g^\prime_\varepsilon(f(z_k,\cdot)-t_k)\nabla f(z_k,\cdot)}{1-\beta}\RQ$. A direct calculation shows that
\begin{equation}\label{eq:variance_without_IS}
\begin{aligned}
\text{Var}_{\rho}\LQ \frac{g^\prime_{\varepsilon}(f(z_k,\cdot)-t_k)\nabla f(z_k,\cdot)}{1-\beta}\RQ &=\E_{\rho}\LQ \left\|\frac{g^\prime_{\varepsilon}(f(z_k,\cdot)-t_k)\nabla f(z_k,\cdot)}{1-\beta} - \E_{\rho}\LQ \frac{g^\prime_{\varepsilon}(f(z_k,\cdot)-t_k)\nabla f(z_k,\cdot)}{1-\beta}\RQ\right\|^2\RQ\\
&=\E_{\rho}\LQ \left\|\frac{g^\prime_{\varepsilon}(f(z_k,\cdot)-t_k)\nabla f(z_k,\cdot)}{1-\beta}\right\|^2 \RQ- \left\|\E_{\rho}\LQ \frac{g^\prime_{\varepsilon}(f(z_k,\cdot)-t_k)\nabla f(z_k,\cdot)}{1-\beta}\RQ\right\|^2\\
&\approx \frac{1}{1-\beta}\E_{\rho}\LQ \left\| \nabla f(z_k,\cdot)\right\|^2|G_{z_k,t_k}\RQ-
\left\|\E_{\rho} \LQ \nabla f(z_k,\cdot)|G_{z_k,t_k}\RQ\right\|^2.
\end{aligned}
\end{equation}
Hence, $\text{Var}_{\rho}\LQ \frac{g^\prime_{\varepsilon}(f(z_k,\cdot)-t_k)\nabla f(z_k,\cdot)}{1-\beta}\RQ$ does not correspond to the variance of the gradient $\nabla f(z_k,\cdot)$ conditioned on the event $G_{z_k,t_k}$, as one may have intuitively guessed since the sampling is not restricted neither to $G_{z_k,t_k}$ or to its smoothed version $G^{\varepsilon}_{z_k,t_k}$.
On the contrary, as $\beta\rightarrow 1$, $\text{Var}_{\rho}\LQ \frac{g^\prime_{\varepsilon}(f(z_k,\cdot)-t_k)\nabla f(z_k,\cdot)}{1-\beta}\RQ$ is dominated by $\frac{1}{1-\beta}\E_{\rho}\LQ \left\| \nabla f(z_k,\cdot)\right\|^2|G_{z_k,t_k}\RQ$, i.e., the second moment of the gradient conditioned to the risk-region. As we will show in Sec. \ref{sec:numerical},  
$\text{Var}_{\rho}\LQ \frac{g^\prime_{\varepsilon}(f(z_k,\cdot)-t_k)\nabla f(z_k,\cdot)}{1-\beta}\RQ$ can be orders of magnitude larger than $\text{Var}_{\rho}\LQ \nabla f(z_k,\cdot)\RQ$, thus leading to much larger values of $M_k$ compared to the risk-neutral setting.

In the next section, we propose to dynamically construct biasing densities $\rhow_k$ from which sampling the set $S_k$, with the goal of reducing the right hand side of \eqref{eq:variance_without_IS}, and thus use smaller sets $S_k$ at each iteration. As a matter of fact, we will show that the variance of the gradients with respect to the adapted biasing distribution satisfies
\begin{equation}\label{eq:claim}
\text{Var}_{\rhow_k}\LQ \frac{g^\prime_{\varepsilon}(f(z_k,\cdot)-t_k)\nabla f(z_k,\cdot)}{1-\beta}\RQ \approx \text{Var}_{\rho}\LQ \nabla f(z_k,\cdot)| G_{z_k,t_k}\RQ,
\end{equation}
These biasing densities will be dynamically constructed along the optimization
and will be adapted to the current risk-region $G^{\varepsilon}_{z_k,t_k}$.

\section{An adaptive importance sampling algorithm}\label{sec:adaptive_importance_sampling}
In this section we present our adaptive algorithm that leverages importance sampling and reduced order models to largely reduce the computational cost of Alg. \ref{alg:adaptive_sampling}.

\subsection{Importance Sampling}
Importance Sampling (IS) is a classical variance-reduction technique to accelerate the estimation of expectations through Monte Carlo methods, see, e.g., \cite[Chapter 5]{asmussen2007stochastic}. The driving idea is to consider a biasing distribution $\rhow$, satisfying $\supp(X\rho)\subset \supp(\rhow)$, which should increase the sampling frequency of events that contribute the most to the estimation $\E\LQ X\RQ$. Then, the integral identity
\[E\LQ X\RQ =\int_{\Gamma} X(\xib)\rho(\xib)d\xib = \int_{\Gamma} X(\xib)w(\xib)\rhow(\xib) d\xib, \]
where $w(\xib):=\frac{\rho(\xib)}{\rhow(\xib)}$, suggests to approximate $\E\LQ X\RQ$ using the \textit{unbiased} Monte Carlo estimator $\frac{1}{M} \sum_{j=1}^M X(\xib_j)w(\xib_j)$, where $\left\{\xib_j\right\}_{j=1}^M$ are i.i.d. samples from the $\rhow$ density, and $\left\{w(\xib_j)\right\}_{j=1}^M$ are weights compensating the change in the sampling distribution. By choosing suitably $\rhow$, one can reduce the variance of the Monte Carlo estimator and obtain more accurate estimates at a reduced number of samples.

In our context, we would like to leverage IS to reduce the variance of the empirical gradient $\nabla \Ftilde^{\varepsilon}_{S_k}$ at every iteration of Alg.\ref{alg:adaptive_sampling}.  Since we have remarked that for every couple $(z_k,t_k)$

\begin{equation*}
\supp \left(\frac{g^{\prime}_{\varepsilon}(f(z_k,\xib)-t_k)\nabla f(z_k,\xib)\rho(\xib) }{1-\beta}\right)\subset G^\varepsilon_{z_k,t_k}\LQ f\RQ,    
\end{equation*}
a reasonable choice for the biasing density would be
\begin{equation*}
\rhow(\xib)=\frac{\mathbbm{1}_{G^{\varepsilon}_{z_k,t_k}\LQ f\RQ}(\xib)\rho(\xib)}{\PP(G^{\varepsilon}_{z_k,t_k}\LQ f\RQ)},    
\end{equation*}
which precisely attributes a non-zero, re-weighted, probability only to samples in $G^{\varepsilon}_{z_k,t_k}\LQ f\RQ$. For a given set of i.i.d. samples $S=\left\{\xib_j\right\}_{j=1}^M$ drawn from $\rhow$, the IS estimator for $\nabla \Ftilde^{\varepsilon}(z_k,t_k)$ would then be
\begin{equation*}
\nabla \Ftilde^{\varepsilon,\ISt}_{S_k}(z_k,t_k)=\frac{1}{M}\sum_{j=1}^M \frac{g^\prime_{\varepsilon}(f(z_k,\xib_j)-t_k)\nabla f(z_k,\xib_j)w(\xib_j)}{1-\beta},
\end{equation*}
where $w(\xib)=\PP(G^{\varepsilon}_{z_k,t_k}\LQ f\RQ)$. However, this promising approach suffers from few drawbacks. First of all, the generation of samples from $\rhow$ may be cumbersome. A standard acceptance-rejection algorithm generates a sample $\xib$ from, e.g., the $\rho$ density, evaluates $f(z_k,\xib)$ and accepts $\xib$ if and only if $\xib\in G^{\varepsilon}_{z_k,t_k}\LQ f\RQ$. This approach can be very computationally expensive whenever the evaluation of $f$ is costly, and it even worsens as $\beta$ tends to one as this leads to higher rejection rates. Second, we will be interested in sampling the set $S_{k}$ when only $z_{k}$ is available, while $t_{k}$ still has to be computed (as it is the case in Alg. \ref{alg:adaptive_sampling}). Of course, since $t_{k}$ is a sampled-based approximation of $\Var \LQ f(z_k,\cdot)\RQ$, we could estimate it through several evaluations of $f(z,\cdot)$, which however may be again computationally expensive, and whose cost can actually be avoided as will show.

To overcome these computational challenges, we leverage Reduced Order Models (ROMs) to both estimate $t_{k+1}$ and then check whether $\xib$ belongs to a suitable risk-region with a reduced computational cost. This in turn induces a ROM-based biasing distribution. We do update the ROM to keep track of the evolution of the control variable $z$ along the optimization, so that a sequence of biasing densities $\left\{\rhow_k\right\}_{k\in \mathbb{N}}$ is constructed.
In the next paragraph, we recall the framework introduced in \cite{heinkenschloss2018conditional} to estimate $\Cvar$ using ROMs, and extend it to the optimization context.

\subsection{Construction of the ROM-based biasing density}\label{sec:ROM_based_density}
We assume that for every $z\in \C$, we have available a Reduced Order Model (ROM) of the map $f(z,\cdot):\Gamma \rightarrow \setR$ denoted by $f_r(z,\cdot)$,  which is much faster to evaluate than $f(z,\cdot)$, and of a computable a-posteriori error bound $e_z:\Gamma \rightarrow \setR^+$ such that
\begin{equation}\label{eq:error_bound_ROM}
|f(z,\xib)-f_r(z,\xib)|\leq e_z(\xib),\quad \rho\text{-a.e. }\xib \in \Gamma.
\end{equation}
Suppose now that we are at the iteration $k$, we have computed $z_{k+1}$ and for the moment, assume further that we can compute $\Var\LQ f(z_{k+1},\cdot)\RQ$ and $\Var\LQ f_r(z_{k+1},\cdot)\RQ$ exactly, and denote them by $\that_{k+1}$ and $\that^r_{k+1}$, respectively. 

A first naive approach to construct a ROM-based biasing distribution from which to sample the next set $S_{k+1}$ would be to consider the risk-region
\[G^{\varepsilon}_{z_{k+1},\that^r_{k+1}}\LQ f_r \RQ :=\left\{\xib: f_r(z_{k+1},\xib)\geq \that^r_{k+1}-\frac{\varepsilon}{2}\right\},\quad\]
and consequently to define
\begin{equation}\label{eq:biasing_density_ROM_wrong}
\rhow(\xib)=\frac{\mathbbm{1}_{G^{\varepsilon}_{z_{k+1},\that^r_{k+1}}\LQ f_r\RQ}(\xib)\rho(\xib)}{\PP(G^{\varepsilon}_{z_{k+1},\that^r_{k+1}}\LQ f_r\RQ)}.  \end{equation}
However, $G^{\varepsilon}_{z_{k+1},\that_{k+1}}\LQ f\RQ \subset G^{\varepsilon}_{z_{k+1},\that^r_{k+1}}\LQ f_r \RQ $ does not necessarily hold, and thus \eqref{eq:biasing_density_ROM_wrong} is not an admissible biasing density. 
Driven by \cite{heinkenschloss2018conditional}, we instead define $\that^{r,e}_{k+1}:=\Var\LQ f_r(z_{k+1},\cdot)-e_{z_{k+1}}\RQ$,
\begin{equation*}
G^{\varepsilon,e}_{z_{k+1},\that_{k+1}^{r,e}}\LQ f\RQ := \left\{\xib:\; g_{\varepsilon}^\prime(f_r(z_{k+1},\xib)+e_{z_{k+1}}(\xib)-\that_{k+1}^{r,e})>0 \right \} = \left\{\xib:\; f_r(z_{k+1},\xib)+e_{z_{k+1}}(\xib)\geq \that_{k+1}^{r,e}-\frac{\varepsilon}{2} \right\},  \end{equation*}
and set
\begin{equation}\label{eq:biasing_density_ROM_right}
\rhow(\xib):=\frac{\mathbbm{1}_{G^{\varepsilon,e}_{z_{k+1},\that^{r,e}_{k+1}}\LQ f_r\RQ}(\xib)\rho(\xib)}{\PP(G^{\varepsilon,e}_{z_{k+1},\that^{r,e}_{k+1}}\LQ f_r\RQ)}.  \end{equation}
The next proposition is adapted from \cite[Lemma 3.2]{heinkenschloss2018conditional}.
\begin{prop}\label{thm:inclusion_risk_regions}
The following inclusion holds
\[G^\varepsilon_{z_{k+1},\that_{k+1}}\LQ f\RQ\subset  G^{\varepsilon,e}_{x_{k+1},\that_{k+1}^{r,e}}\LQ f\RQ.\] 
\end{prop}
\begin{pf}
From \eqref{eq:error_bound_ROM} it follows that $f(z_{k+1})\geq f_r(z_{k+1},\cdot)-e_{z_{k+1}}$ which implies $\that_{k+1}=\Var\LQ f(z_{k+1},\cdot)\RQ\geq \Var\LQ f_r(z_{k+1},\cdot)-e_{z_{k+1}}\RQ =\that^{r,e}_{k+1}$ due to the monotonicity of $\Var $. Then, taking a $\xib\in G^\varepsilon_{z_{k+1},\that_{k+1}}\LQ f\RQ$, the sequence of inequalities
\[f_r(z_{k+1},\xib)+e_{z_{k+1}}(\xib)\geq f(z_{k+1},\xib)\geq \that_{k+1}-\frac{\varepsilon}{2}\geq \that^{r,e}_{k+1} -\frac{\varepsilon}{2},\] 
shows that $\xib\in G^{\varepsilon,e}_{z_{k+1},\that_{k+1}^{r,e}}\LQ f_r\RQ$.
\qed 
\end{pf}
Proposition \ref{thm:inclusion_risk_regions} guarantees that \eqref{eq:biasing_density_ROM_right} is a suitable biasing density. However, both $\that_{k+1}$ and $\that_{k+1}^{r,e}$ are not computable quantities. In a numerical implementation, at iteration $k$ we underestimate the future $t_{k+1}$ (which is computed by minimizing $\Ftilde^{\varepsilon}(z_{k+1},t)$ with respect to $t$ at iteration $k+1$) via
\begin{equation*}
t^{r,e}_{k+1}=\argmin_{t\in\setR}\; t+\frac{1}{(1-\beta)M_{\text{Trial}}}\sum_{j=1}^{M_{\text{Trial}}} g^{\varepsilon}\left(f_r(z_{k+1},\xib^{\text{Trial}}_j)-e_{z_{k+1}}(\xib^{\text{Trial}}_j)-t\right),
\end{equation*}
where $\left\{\xib^{\text{Trial}}_j\right\}_{j=1}^{M_{\text{Trial}}}$ are distributed according to $\rho$. The new i.i.d. samples $S_{k+1}$ are then drawn from the biasing density 
\begin{equation*}
\rhow^{k+1}(\xib):=\frac{\mathbbm{1}_{G^{\varepsilon,e}_{z_{k+1},t^{r,e}_{k+1}}\LQ f_r\RQ}(\xib)\rho(\xib)}{\PP(G^{\varepsilon,e}_{z_{k+1},t^{r,e}_{k+1}}\LQ f_r\RQ)}.
\end{equation*}
In principle, it might occur that the estimated $t^{r,e}_{k+1}$ is larger than the actual $t_{k+1}$ (computed at iteration $k+1$). In our numerical simulations we never encountered such case which, however, can be easily avoided by recomputing $t^{r,e}_{k+1}$ with a larger value of $M_{\text{Trial}}$ (as asymptotically $t^{r,e}_{k+1}$ converges to $\that^{r,e}_{k+1}$, hence it will be eventually smaller than $t_{k+1}$ \cite{hong2014monte}), followed by another acceptance-rejection algorithm to generate a new set $S_{k+1}$.

\subsection{The algorithm}
Our adaptive importance sampling algorithm is described in Alg. \ref{alg:adaptive_sampling_rom}. At every iteration, an alternating procedure (Steps 5-6) identical to Alg. \ref{alg:adaptive_sampling} is carried out, although now the IS estimator 
\begin{equation*}
\Ftilde^{\varepsilon,\ISt}_{S_k}(z,t):=t+\frac{1}{(1-\beta)M_k}\sum_{\xib_k^j\in S_k} g_{\varepsilon}(f(z,\xib_j^k)-t)w_k,
\end{equation*} is employed, with weights $w_k=\PP\left( G^{\varepsilon,e}_{z_{k+1},t^{r,\epsilon}_{k+1}}\LQ f_r\RQ\right)$. Condition \eqref{eq:norm_test_the} now reads
\begin{equation*}
\widetilde{\E}_k\LQ \|\nabla \Ftilde^{\varepsilon,\ISt}_{S_k}(z_k,t_k)-\nabla \Ftilde^{\varepsilon}(z_k,t_k)\|^2\RQ =\frac{\Variance_{\rhow_k}\LQ\frac{g_\varepsilon^\prime(f(z_k,\cdot)-t_k)\nabla f(z_k,\cdot)w_k}{1-\beta}\RQ}{M_k} \leq \theta^2 \|\nabla \Ftilde^{\varepsilon}(z_k,t_k)\|^2,
\end{equation*}
where $\widetilde{\E}_k$ is the expectation with respect to the set of $M_k$ i.i.d. random variables distributed according to $\rhow_k$. The practical norm test is thus,
\begin{equation}\label{eq:norm_test_cvar_pra_IS}
\varrho=\frac{\sum_{\xib^k_j\in S_k}\left\| \frac{g^\prime_{\varepsilon}(f(z_k,\xib_j^k)-t_k)\nabla f(z_k,\xib_j^k)w_k}{1-\beta}-\nabla \Ftilde^{\varepsilon,IS}_{S_k}(z_k,t_k)\right\|^2}{(M_k-1)(M_k)\theta^2 \|R_{S_k}(z_k)\|^2}.
\end{equation}
Notice that a direct calculation shows
\begin{equation}\label{eq:variance_with_IS}
\Variance_{\rhow_k}\LQ\frac{g_\varepsilon^\prime(f(z_k,\cdot)-t_k)\nabla f(z_k,\cdot)w_k}{1-\beta}\RQ\approx \frac{w_k}{1-\beta}\E_{\rho}\LQ \left\| \nabla f(z_k,\cdot)\right\|^2|G^{\varepsilon}_{z_k,t_k}\RQ-
\left\|\E_{\rho} \LQ \nabla f(z_k,\cdot)|G_{z_k,t_k}\RQ\right\|^2,
\end{equation}
where the approximation is due to the fact that $g^\prime_{\varepsilon}$
is a smoothed version of $\mathbbm{1}_{G^{\varepsilon}_{z_k,t_k}}$ and that $t_k$ is only a sample-based estimate of $\Var\LQ f(z_k,\cdot)\RQ$.
It follows that, if along the optimization the ROM model becomes increasingly accurate, then $w_k$ tends to $(1-\beta)$, and we obtain claim \eqref{eq:claim}. We refer to Section \ref{sec:numerical} for a quantitative comparison between \eqref{eq:variance_without_IS} and \eqref{eq:variance_with_IS}.

Once the next sample size has been computed, we construct a reduced order model for the new iterate $z_{k+1}$. We then draw a set of $M_{\text{Trial}}$ samples, evaluate the a-posteriori error estimator for each of these samples, compute $t^{r,e}_{k+1}$ and estimate the size of the region $G^{\varepsilon,e}_{z_{k+1},t_{k+1}^{r,e}}$ (lines 11-14). Finally, a standard acceptance-rejection algorithm is used to generate samples $\left\{\xib_j\right\}_{j=1}^{M_{k+1}}$ distributed according to $\rhow^{k+1}$. Notice that both the computation of $t^{r,e}_{k+1}$ and the test if $\xib_{\text{test}}$ belongs to $G^{\varepsilon,e}_{z_{k+1},t_{k+1}^{r,e}}$ rely exclusively on the ROM $f_r$, so that the full-order model is never used between lines 11 and 23.
\begin{algorithm}
\caption{Adaptive Importance Sampling}\label{alg:adaptive_sampling_rom}
\begin{algorithmic}[1]
\State \textbf{Require} $\mathbf{u}^0$, $M_0$, $M_{\max}$, $\alpha$, $M_{\text{trial}}$
\State Set $k=0$.
\State Sample $S_k=\left\{\xib_i\right\}_{i=1}^{M_k}$ from $\rho(\xib)$.
\While {$\sum_{j=0}^k M_j\leq M_{\max}$}
\State Compute $t_{k}\in \text{argmin}_{t\in \mathbb{R}} F_{S_k}^{\varepsilon,\ISt}(z_k,t)$.
\State $z_{k+1}=z_k-\alpha \nabla F^{\varepsilon,\ISt}_{S_k}(z_k,t_{k})$.
\State Compute $\varrho$ in \eqref{eq:norm_test_cvar_pra_IS}.
\If{$\varrho>1$} $M_{k+1}=\varrho M_k$ \EndIf.
\State Construct ROM $f_r(z_{k+1},\cdot)$
\State Draw $M_{\text{Mtrial}}$ samples from $\rho$: $S_{\text{Trial}}=\left\{\xib_j\right\}_{j=1}^{M_{\text{trial}}}$.
\State Compute a-posteriori error estimator: $|f(z_{k+1},\xib_j)-f_r(z_{k+1},\xib_j)|\leq e_{z_{k+1}}(\xib_j)$, $j=1,\dots,M_{\text{Trial}}$.
\State Compute $t^{r,e}_{k+1}=\argmin_{t\in\setR} t+\frac{1}{(1-\beta)M_{\text{Trial}}}\sum_{j=1}^{M_{\text{Trial}}} g^{\varepsilon}\left(f_r(z_{k+1},\xib_j)-e_{z_{k+1}}(\xib_j)-t\right)$.
\State Estimate  $\PP\left(G^{\varepsilon,e}_{z_{k+1},t^{r,e}_{k+1}}\LQ f_r\RQ\right)=\#\left\{\xib\in S_{\text{Trial}}:\xib \in G^{\varepsilon,e}_{z_{k+1},t^{r,e}_{k+1}}\right\}/M_{\text{Trial}}$.
\State Importance sampling step: $n=0$.
\While {$n\leq M_{k+1}$}
\State Draw a sample $\xib_{test}$ from density $\rho$.
\If{$\xib_{test}\in G^{\varepsilon,e}_{z_{k+1},t^{r,e}_{k+1}}\LQ f_r\RQ$} 
\State Add sample $\xib_{test}$ to $S_{k+1}$. 
\State $n\gets n+1$. 
\EndIf
\EndWhile
\EndWhile
\end{algorithmic}
\end{algorithm}
We conclude this section remarking that the efficient construction of a ROM along the optimization is highly problem dependent and may be nontrivial. 
In Section \ref{sec:numerical}, we discuss in detail an approach in the context of PDE-constrained optimization under uncertainty. Here we limit to make the interesting, and maybe surprising, observation: the ROM does not need to be particularly over accurate in the risk-region. While surely an overall accurate ROM leads to sharper estimates of $t_{k+1}^{r,e}$, the ROM is used here essentially to exclude realization of $\xib$ that do not contribute to the $\Cvar$. Thus, large values of the a-posteriori error estimator outside the ``true" risk region $G^{\varepsilon}_{z_k,t_k}$ are those problematic, since they may actually lead to the acceptance of unnecessary samples in the acceptance-rejection step.

\subsection{Convergence analysis}\label{sec:convergence}
In this subsection we analyze the convergence of Alg. \ref{alg:adaptive_sampling_rom}. Readers more interested in implementation details and numerical results can
proceed directly to Section \ref{sec:numerical}.

We start by remarking that to analyze the convergence of Alg. \ref{alg:adaptive_sampling_rom} it is sufficient to study Alg. \ref{alg:adaptive_sampling}. Indeed, the additional complexity of Alg. \ref{alg:adaptive_sampling_rom} due to the ROM-based construction of the biasing distribution only affects the number of samples needed to satisfy a condition on the gradients (see \eqref{eq:condition_gradients}).
While \cite{beiser2023adaptive} proposes Alg. \ref{alg:adaptive_sampling} without a convergence analysis, a first study of an alternating procedure for the $\Cvar$ minimization has been proposed in \cite{ganesh2023gradient}, under however the hypothesis that $\Ftilde$ is jointly strongly convex for all $(z,t)\in \mathcal{Z}\times \setR$. This assumption is somehow too strong since if $t\notin \text{supp}(f(z,\cdot))$, then $\Ftilde$ is only linear in $t$ in a neighborhood of $(z,t)$. In practice though, Alg. \ref{alg:adaptive_sampling} never explores these problematic regions since the minimization step in $t$ guarantees that $t\in \text{supp}(f(z,\cdot))$.

Our analysis proceeds through a different path, consisting of reinterpreting the alternating procedure of Alg. \ref{alg:adaptive_sampling} as a single gradient step on a reduced functional defined only over the $z$ variable. Under suitable assumptions detailed next, we argue that this reduced functional is strongly convex and has Lipschitz continuous gradients, so that the convergence of an adaptive sampling gradient descent algorithm directly follows from available results in literature.

To begin our analysis, we make the following assumption on the properties of $\Ftilde^{\varepsilon}$ with respect to the variable $t$.
\begin{ass}\label{ass:h}
For every $z\in \C$, the minimization problem 
\[\min\limits_{t\in \setR} \Ftilde^{\epsilon}(z,t)\]
admits an unique minimum $t^\star(z)$ and 
defines a differentiable map $h:\C\rightarrow \setR$ such that $h(z)=t^\star(z)$.
\end{ass}
Notice that, due to the differentiability of $\Ftilde^{\varepsilon}$ with respect to $t$, $\partial_t \Ftilde^{\varepsilon}(z,h(z))=0$ trivially holds and, assuming the uniqueness of $t^\star(z)$, the local existence and smoothness of $h$ follow from the implicit function theorem. Assumption \ref{ass:h} allows us to introduce the reduced functional
\begin{equation*}
\Fhat^{\varepsilon}(z):=h(z)+\frac{1}{1-\beta}\E\LQ g^{\varepsilon}(f(z,\cdot)-h(z))\RQ,    
\end{equation*}
and a straight calculation shows that 
\begin{equation}\label{eq:property_gradients}
\nabla \Fhat^{\varepsilon}(z)=\nabla \Ftilde^{\varepsilon}(z,h(z))+\partial_t \Ftilde^{\varepsilon}(z,h(z))\nabla h(z)=\nabla \Ftilde^{\varepsilon}(z,h(z)),
\end{equation}
where the term $\partial_t \Ftilde^{\varepsilon}$ cancels due to the definition of $h$.
From \eqref{eq:property_gradients}, it follows that
the alternating procedure
\[t_k=\min_{t\in\setR} \Ftilde^{\varepsilon}(z_k,t),\quad z_{k+1}=z_k-\alpha\nabla \Ftilde^{\varepsilon}(z_k,t_k),\]
is equivalent to the iteration
\[\tilde{z}_{k+1}=\tilde{z}_{k}-\alpha \nabla \Fhat^{\varepsilon}(\tilde{z}_k),\]
in the sense that $\tilde{z}_{k+1}=z_{k+1}$ and $t_{k}=h(z_k)$ for every $k$. 
Our goal is then to study the convergence of an adaptive sampling gradient algorithm applied to the reduced functional $\Fhat^\varepsilon$. To do so, we make the following assumptions on $f(z,\cdot)$.
\begin{ass}\label{ass:f}
For $\rho$-a.e. $\xib$, the mapping $z\rightarrow f(z,\xib)$ is twice differentiable, and there exist two positive constants $\mu_f,L_f
\in \setR^+$, such that $\forall \;\delta z\in \mathcal{Z}$
\[ \mu_f\|\delta z\|^2 \leq(\delta z,H_f(z,\xib)\delta z)\leq L_f\|\delta z\|^2\quad \rho\text{-a.e. }\xib,\]
$H_f(z,\xib)$ being the Hessian of $f$ at $(z,\xib)$.
Furthermore, there exists a positive constant $C_f\in \setR$ such that $\forall z\in \C$, $\|\nabla f(z,\cdot)\|\leq C_f$ $\rho$-a.e..
\end{ass}
Assumption \ref{ass:f} implies that $f$ is strongly convex, and has Lipschitz continuous gradients $\rho$-a.s.. This in turn allows us to show that the same properties hold true for $\Fhat^{\varepsilon}$.
\begin{lem}\label{lemma:properties_Fhat}
The functional $\Fhat^{\varepsilon}$ is strongly convex over $\C$ and has Lipschitz continuous gradients.
\end{lem}
\begin{pf}
We first prove the Lipschitz continuity of gradients. Using that  $\sup_{t\in\setR} |g_{\varepsilon}^\prime(t)|\leq 1$, eq. \eqref{eq:property_gradients}, Assumption \ref{ass:f}, and the Lipschitz continuity of $f,h$ and $g^\prime_{\varepsilon}$,
\begin{equation*}
\begin{aligned}
\left\|\nabla \Fhat^{\varepsilon}(x)-\nabla \Fhat^{\varepsilon}(y)\right\|&\leq \frac{1}{1-\beta}\E\LQ \left\| g^\prime_\varepsilon(f(x,\cdot)-h(x))\nabla f(x,\cdot)- g^\prime_\varepsilon(f(y,\cdot)-h(y))\nabla f(y,\cdot)\right\|\RQ\\
&\leq \frac{1}{1-\beta}\E\bigg[ g^\prime_\varepsilon(f(x,\cdot)-h(x))\left\|\nabla f(x,\cdot)-\nabla f(y,\cdot)\right\|\\
&+(g^\prime_\varepsilon(f(x,\cdot)-h(x)) - g^\prime_\varepsilon(f(y,\cdot)-h(y)))\left\|\nabla f(y,\cdot)\right\|\bigg]\\
&\leq \frac{1}{1-\beta}\left(\|x-y\| L_f + C_f\E\LQ g^\prime_\varepsilon(f(x,\cdot)-h(x)) - g^\prime_\varepsilon(f(y,\cdot)-h(y))\RQ\right)\\
&\leq L_{\Fhat^{\varepsilon}}\|x-y\|,
\end{aligned}
\end{equation*}
for a suitable positive constant $L_{\Fhat^{\varepsilon}}$.
Next, we focus on the strong convexity. Denoting with $\Hc_{\Fhat^\varepsilon}$ the Hessian of $\Fhat^\varepsilon$, the proof consists in showing that
\begin{equation*}
\min_{\delta  z\in\mathcal{Z}} \frac{(\delta z, \Hc_{\Fhat^{\varepsilon}}(z) \delta z)}{\|\delta z\|^2}=:\mu_{\Fhat^{\varepsilon}}>0,\quad \forall z\in \mathcal{Z}.
\end{equation*}
To do so, we will actually leverage the Hessian of $\Ftilde^{\varepsilon}$, that is,
\[
\Hc_{\Ftilde^{\varepsilon}}(z,t):=\begin{pmatrix}
\Ftilde^{\varepsilon}_{zz}(z,t) & \Ftilde^{\varepsilon}_{tz}(z,t)\\
\Ftilde^{\varepsilon}_{zt}(z,t) & \Ftilde^{\varepsilon}_{tt}(z,t)
\end{pmatrix}    
\]
where, using the shorthand notation $\bar{g}_{\varepsilon}=g_{\varepsilon}(f(z,\cdot)-h(z))$, direct calculations show that
\begin{equation*}
\begin{aligned}
(\delta z,\Ftilde_{zz}^{\varepsilon}(z,t)\delta z)&=\frac{1}{1-\beta}\E\LQ \bar{g}_{\varepsilon}^{\dprime} (\delta z,\nabla f(z,\cdot))(\nabla f(z,\cdot),\delta z) + \bar{g}_{\varepsilon}^\prime (\delta z, H_f(z,\cdot)\delta z)\RQ,\\
(\delta t, \Ftilde_{zt}^{\varepsilon}(z,t)\delta z)&=(\delta z, \Ftilde_{tz}^{\varepsilon}(z,t) \delta t)=-\frac{\delta t}{1-\beta}\E\LQ \bar{g}_{\varepsilon}^\dprime (\nabla f(z,\cdot),\delta z)\RQ,\\
(\delta t,\Ftilde_{tt}^{\varepsilon}(z,t)\delta t)&=\frac{\delta t^2}{1-\beta}\E\LQ \bar{g}_{\varepsilon}^\dprime \RQ,
\end{aligned}
\end{equation*}
As a matter of fact, it is well-known $\Hc_{\Fhat^\varepsilon}(z)$ corresponds to the Schur complement of $\Hc_{\Ftilde^{\varepsilon}}$ evaluated in $(z,h(z))$, namely
\[
\Hc_{\Fhat^\varepsilon}(z)=\Ftilde^{\varepsilon}_{zz}(z,h(z))-\Ftilde^{\varepsilon}_{tz}(z,h(z))\Ftilde^{\varepsilon}_{tt}(z,h(z))^{-1}\Ftilde^{\varepsilon}_{zt}(z,h(z)),
\]
so that
\begin{equation}\label{eq:step_prof}
\begin{aligned}
(\delta z,\Hc_{\Fhat^\varepsilon}(z) \delta z)&=(\delta z,\Ftilde^{\varepsilon}_{zz}(z,h(z))\delta z)- (\delta z, \Ftilde_{tz}^\varepsilon (z,h(z))\Ftilde_{tt}^{-1}(z,h(z)) \Ftilde_{zt}^\varepsilon(z,h(z)) \delta z)\\
&=\frac{1}{1-\beta} \E\LQ \bar{g}^\prime_\varepsilon(\delta z, H_f(z,\cdot)\delta z)+\bar{g}_\varepsilon^\dprime |(\nabla f(z,\cdot),\delta z)|^2\RQ\\
&-\frac{1}{1-\beta}\frac{\left(\E\LQ \bar{g}_\varepsilon^\dprime (\nabla f(z,\cdot),\delta z)\RQ\right)^2}{\E\LQ \bar{g}_\varepsilon^\dprime\RQ}.
\end{aligned}
\end{equation}
Using the positivity of $\bar{g}^\dprime$ and the Cauchy-Schwarz inequality,
\[\frac{\left(\E\LQ \bar{g}_\varepsilon^\dprime (\nabla f(z,\cdot),\delta z)\RQ\right)^2}{\E\LQ \bar{g}_\varepsilon^\dprime\RQ}\leq \E\LQ \bar{g}_\varepsilon^\dprime |(\nabla f(z,\cdot),\delta z)|^2\RQ\]
which inserted into \eqref{eq:step_prof} leads to
\[(\delta z,\Hc_{\Fhat^\varepsilon}(z)\delta z)\geq \frac{1}{1-\beta} \E\LQ \bar{g}^\prime_\varepsilon(\delta z, H_f(z,\cdot)\delta z)\RQ \geq \frac{\mu_f}{1-\beta}\E\LQ \bar{g}^\prime_{\varepsilon}(f(z,\cdot)-h(z))\RQ\|\delta z\|^2\geq \frac{\mu_f}{2}\|\delta z\|^2,\]
where the last step uses that $\E\LQ \bar{g}^\prime_{\varepsilon}(f(z,\cdot)-h(z))\RQ\geq \frac{1}{2}\E\LQ \mathbbm{1}_{f(z,\cdot)-h(z)}\RQ=\frac{1-\beta}{2}.$
\qed 
\end{pf}

We are now almost ready to state a convergence result for Alg. \ref{alg:adaptive_sampling}. We however need an additional assumption on the sample function $F_{S_k}$.
\begin{ass}\label{ass:h_sampled}
We assume that for every $z\in \C$ the minimization problem 
\[\min\limits_{t\in \setR} F_{S_k}^{\epsilon}(z,t)\]
admits a unique minimum $t_{S_k}^\star(z)$ and 
defines a map $h_{S_k}:\C\rightarrow \setR$ such that $h_{S_k}(z)=t_{S_k}^\star(z)$.
\end{ass}
Assumption \ref{ass:h_sampled} allows us to introduce the reduced version of $F_{S_k}$,
\[\Fhat^\varepsilon_{S_k}(z):=h_{S_k}(z)+\frac{1}{1-\beta}\E\LQ g_{\varepsilon}(f(z,\cdot)-h_{S_k}(z))\RQ.\]
Repeating verbatim the argument made for $\Fhat^\varepsilon$, the alternating procedure of Alg. \ref{alg:adaptive_sampling} is equivalent to the gradient descent iteration
\begin{equation}\label{eq:iteration_sample_approximation}
z_{k+1}=z_k-\alpha \nabla \Fhat^\varepsilon_{S_k}(z_k).    
\end{equation}
The next theorem shows that choosing at each iteration a sample $S_k$ such that
\begin{equation}\label{eq:condition_gradients}
\E_{k}\LQ \left\|\nabla \Fhat^{\varepsilon}_{S_k}(z_k)-\nabla \Fhat^{\varepsilon}(z_k)\right\|^2\RQ\leq \theta^2\|R(z_k)\|^2
\end{equation}
is sufficient to preserve linear convergence of the iterates $\left\{z_k\right\}_{k\geq 1}$.
\begin{thm}\label{thm:convergence}
Let $\left\{z_k\right\}_{k\geq 1}$ be the sequence generated by \eqref{eq:iteration_sample_approximation}, with $\alpha< L_{\Fhat^\varepsilon}$ and such that each $S_k$ satisfies \eqref{eq:condition_gradients}. Then, for sufficiently small $\theta>0$, $z_k$ converges linearly in expectation, i.e.,
\[
\E\LQ \|z_{k+1}-z^\star\|\RQ \leq \varrho^k\|z_0-z^\star\|,    
\]
for some $\varrho\in [0,1).$
\end{thm}
\begin{pf}
The claim follows directly from \cite[Theorem 2.10]{beiser2023adaptive} since due to Lemma \ref{lemma:properties_Fhat}, $\Fhat^\varepsilon$ is strongly convex and has Lipschitz continuous gradients.  
\qed 
\end{pf}

\begin{rmk}
Similarly to \eqref{eq:norm_test_the}, Condition \eqref{eq:condition_gradients} is not numerically implementable since $\Fhat$ is not computable. However, in addition $\nabla \Fhat^{\varepsilon}_{S_k}(z_k)$ is a biased estimator of $\nabla \Fhat^{\varepsilon}(z_k)$, since $h_{S_k}$ is only an asymptotic unbiased estimator of $h$, see, e.g., \cite{hong2014monte}. Thus,
\[\E_{k}\LQ \left\|\nabla \Fhat^{\varepsilon}_{S_k}(z_k)-\nabla \Fhat^{\varepsilon}(z_k)\right\|^2\RQ=\E_k\LQ \left\|\nabla \Fhat^{\varepsilon}_{S_k}(z_k)-\E_k\LQ \nabla \Fhat^{\varepsilon}_{S_k}(z_k)\RQ\right\|^2\RQ + \E\LQ \left\|\E\LQ \nabla \Fhat^{\varepsilon}_{S_k}(z_k)\RQ -\nabla \Fhat^{\varepsilon}(z_k)\right\|^2\RQ,\]
where the first term corresponds to the variance of $\nabla \Fhat_{S_k}(z_k,t_k)$, while the second contribution is a bias term that could be further bounded using the Lipschitz continuity of $g_\varepsilon^\prime$ and Assumption \ref{ass:f}, 
\[\E_k\LQ \left\|\E\LQ \nabla \Fhat^{\varepsilon}_{S_k}(z_k)\RQ -\nabla \Fhat^{\varepsilon}(z_k)\right\|^2\RQ\leq
\frac{C_f^2}{(1-\beta)^2}\E_k\LQ|h(z_k)-h_{S_k}(z_k)|^2\RQ.\]
In our numerical experiments, we neglect the bias term and choose the next sample size by exclusively control the variance term using the practical norm test \eqref{eq:norm_test_cvar_pra_IS}.
\end{rmk}

\section{Numerical Experiments}\label{sec:numerical}
In this section we present numerical experiments arising from the field of PDE-constrained optimization under uncertainty. 
We consider the domain $\D=(0,1)^2$ whose boundary $\partial \Omega$ is divided into two disjoint sets, $\Gamma_D:=(0,1)\times \left\{0,1\right\}$ and $\Gamma_N:=\partial \Omega\setminus \Gamma_D$, on which homogeneous Dirichlet and Neumann boundary conditions are imposed. We denote by $V$ the Sobolev space $H^1_{0,\Gamma_D}(\D)$. The PDE-constrained optimization problem reads
\begin{equation}\label{eq:model_problem_num}
\min_{z\in \mathcal{Z}} \Cvar\LQ \frac{1}{2}\| y(z,\xib)\|^2_{L^2(\D_0)}+\frac{\nu}{2}\|z\|^2_{\mathcal{Z}}\RQ,    
\end{equation}
where $\D_0:=(0.5,1)\times (0,1) \subsetneq \D$ and $y(z,\xib)$ is the solution of the random partial differential equation,
\begin{equation}\label{eq:state_equation}
\int_{\D}\kappa(\xb,\xib) \nabla y(\xb;z,\xib) \cdot\nabla v(\xb)\;d\xb=\int_{\D} l(\xb)v(\xb)\;d\xb+\langle Bz,v\rangle,\quad \forall v\in V.    
\end{equation}
Here, $l(\xb)$ models three point source terms,
\[l(\xb):=10^2\sum_{i=1}^3 \exp\left(-\frac{\|\xb -\xb_i\|}{L}\right),\]
where $L=0.05$, $\xb_j=[0.25,j 0.25]$ for $j=1,\dots,3$, and $B:\mathcal{Z}\rightarrow V^\prime$ is a control operator representing a Neumann boundary control acting on the right vertical edge (denoted from now $E$),
\[\langle Bz, v\rangle=\int_E z(\xb)v(\xb)\;d\xb,\]
hence $\mathcal{Z}=L^2(E)$. Problem \eqref{eq:model_problem_num} is a simple model of a heating process where the heat generated at specific locations has to be extracted on the right portion of the domain by applying a suitable flux along the right edge. It can be embedded into the notation of the previous sections by defining the map $f:\mathcal{Z}\times (0,1)^2\rightarrow \setR$, with $f(z,\xib)=\frac{1}{2}\|y(z,\xib)\|^2_{L^2(\Omega_0)}+\frac{\nu}{2}\|z\|^2_{L^2(E)}$.
 
\begin{figure}
    \centering
    \includegraphics[width=0.3\linewidth]{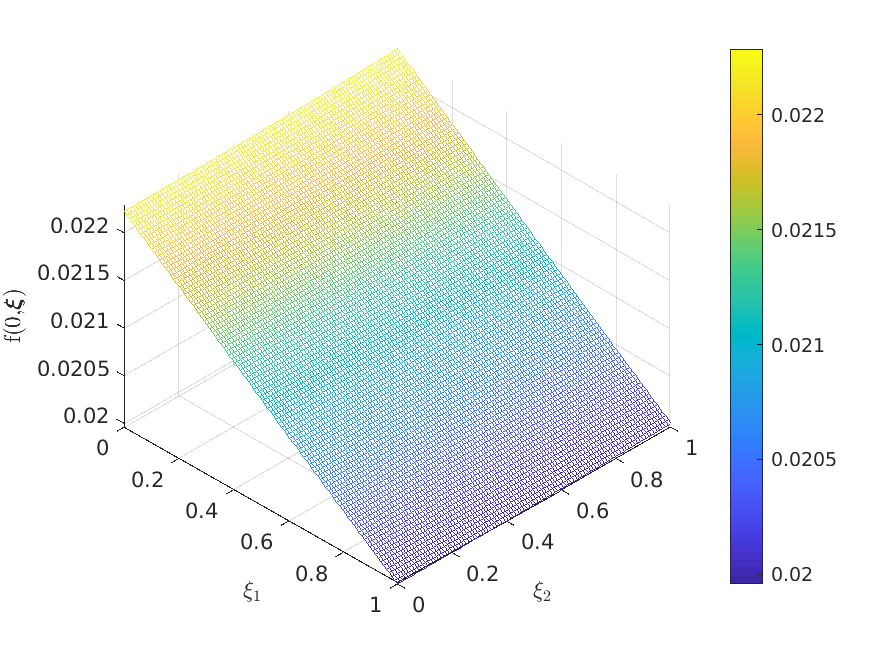}
    \includegraphics[width=0.3\linewidth]{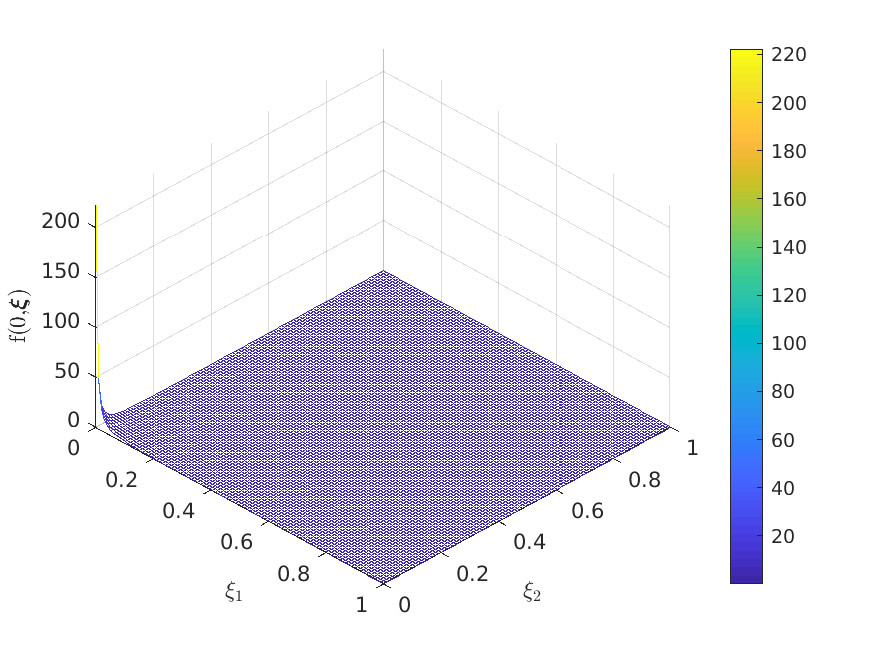}
    \caption{Graphical representation of the map $\xib\rightarrow f(0,\xib)$ for $\kappa_1$ and $\kappa_2$.}
    \label{Fig.Qoi}
\end{figure}

We will consider two different diffusion coefficients, namely
\[\kappa_1(\xb,\xib):= 1+\frac{\xi_1}{2}\sin(\pi x)\sin(\pi y)+\frac{\xi_2}{5}\sin(\pi x)\sin(2\pi y),\]
representing a truncated Fourier expansion, and
\[\kappa_2(\xb,\xib):= 0.01+\xi_1 \exp\left(\sin(\pi x)\sin(\pi y)\right)+\xi_2 \exp\left(\cos(\pi x)\cos(\pi y)\right).\]
The parameters $\xib=(\xi_1,\xi_2)$ are uniformly distributed in $(0,1)^2$. We choose two random variables $(M=2)$ so that an accurate reference solution of $\eqref{eq:model_problem_num}$ can be computed using a tensorized Gauss-Legendre quadrature formula \cite{babuvska2007stochastic}. We emphasize that the two diffusion coefficients lead to very different properties of the objective function. Figure \ref{Fig.Qoi} shows the values of the map $\xib\rightarrow f(0,\xib)$, that is of the uncontrolled ($z=0$) physical process. For $\kappa_1$, $f(0,\cdot)$ has a favorable behavior with overall minor variations across the parameters domain, whereas for $\kappa_2$, $f(0,\cdot)$ assumes very large values within a very tiny region around the vertex $(0,0)$. This is an extremely challenging case for adaptive sampling algorithms since the uniform random sampling rarely generates samples very close to $(0,0)$, which however are the dominant contributions to the $\Cvar$ and to the variance of the sampled gradients.

\subsection{Dynamic update of the ROM}
In this subsection we detail how the ROMs $f_r(z_{k+1},\cdot)$ are computed on the fly during Alg. \ref{alg:adaptive_sampling_rom} (line 10) at a negligible computational cost.

Let $\left\{V_r^{k+1}\right\}_{k\geq 0}$ be a sequence of finite dimensional subspaces, not necessarily nested, of $V$. At a given iterate $(z_{k+1},t_k)$, we consider the ROM
\begin{equation*}
    f_r(z_{k+1},\cdot):=\frac{1}{2}\|y_r(z_{k+1},\cdot)\|^2_{L^2(\Omega_0)}+\frac{\nu}{2}\|z_{k+1}\|^2_{\mathcal{Z}},
\end{equation*}
where $y_r(z_{k+1},\xib)\in V_r^{k+1}$ is the solution of the reduced basis equation \cite{quarteroni2015reduced,hesthaven2016certified},
\[\int_{\D}\kappa(\xb,\xib) \nabla y_r(\xb;z_{k+1},\xib) \cdot\nabla v(\xb)\;d\xb=\int_{\D} l(\xb)v(\xb)\;d\xb+\langle Bz_{k+1},v\rangle,\quad \forall v\in V_r^{k+1}.\]
Clearly, this approach is computationally sound only if the subspace $V_r^{k+1}$ is of small dimension and it can be obtained cheaply. To achieve this, we leverage that the calculation of $t_k$ requires the evaluation of $f(z_k,\cdot)$ at all realizations $\xib\in S_k$. This in turn implies that a solution $y(z_k,\xib)$ to \eqref{eq:state_equation} is computed for each $\xib\in S_k$. We therefore consider a subset $\I_1\subset S_k$ of the snapshots $\left\{y(z_k,\xib)\right\}_{\xib\in S_k}$, and extract the most significant modes using Proper Orthogonal Decomposition, see, e.g., \cite[Section 6]{quarteroni2015reduced}. 
However, as discussed in Section \ref{sec:ROM_based_density}, the ROM needs to be accurate also in the complement of the risk-region. Hence, we enhance the POD basis using a Greedy approach that 
relies on a set of samples $\I_2$ that lie outside the current risk-region $G_{z_{k+1},t^{r,e}_{k+1}}\LQ f_r\RQ$. The set $\I_2$ can be readily obtained by saving a few of the rejected samples during the acceptance-rejection construction of the biasing density.
The Greedy procedure consists in iteratively adding the snapshots that maximize a reduced basis a-posterior error estimator $e^{\RB}_{z_{k+1}}$, which measures how well the current subspace $V_r^{k+1}$ approximates the solution manifold of \eqref{eq:state_equation}. In our implementation, we used a classical residual-based error estimator, see, e.g., \cite[Chapter 3]{quarteroni2015reduced}. The overall procedure to construct a ROM at each iteration is described by Alg. \ref{alg:greedy}.
In our numerical implementation, we provide to Alg. \ref{alg:greedy} two sets $\I_1$ and $\I_2$ consisting of, respectively, $\min\left\{M_k,40\right\}$ randomly chosen elements of $S_k$, and of maximum 100 samples that lie outside $G_{z_{k+1},t^{r,e}_{k+1}}\LQ f_r\RQ$. The tolerance is decreased along the outer iterations of Alg. \ref{alg:adaptive_sampling_rom} as we wish to have a more accurate ROM as we approach convergence.
\begin{algorithm}
\caption{Greedy-enhanced construction of $f_r(z_{k+1},\cdot)$}\label{alg:greedy}
\begin{algorithmic}[1]
\State \textbf{Require} $\left\{y(z_k,\xib)\right\}_{\xib\in \I_1}$, $\I_2$, $\Tol\in\setR^+$, $e^{\RB}_{z_{k+1}}$.
\State $V^{k+1}_r$=$\text{POD}(\left\{y(z_k,\xib)\right\}_{\xib\in \I_1},\Tol)$.
\While {$\max_{\xib\in \I_2}e^{\RB}_{z_{k+1}}(\xib)\geq \Tol$ }
\State Choose $\xib^\star=\argmin_{\xib\in \I_2}e^{\RB}_r{z_{k+1}}(\xib)$.
\State $V_r^{k+1}=V_r^{k+1}\cup y(z_{k+1},\xib^\star)$.
\State Set $\I_2=\I_2\setminus \xib^\star$.
\EndWhile
\State Return $V_r^{k+1}$.
\end{algorithmic}
\end{algorithm}

Notice that, despite $V_r^{k+1}$ is constructed using snapshots corresponding to both the iterates $z_k$ and $z_{k+1}$, a standard residual-based a-posteriori error still provides a rigorous upper bound for the error $\|y(z_{k+1},\cdot)-y_r(z_{k+1},\cdot)\|^2_{L^2(\D)}$. However, the latter is generally non zero even at realizations whose snapshots are included in $V_r^{k+1}$. 
The availability of $e^{\RB}_{z_{k+1}}$ permits to 
estimate the a-posterior error estimator $e_{z_{k+1}}$ for the ROM model since
\begin{equation*}
\begin{aligned}
|f(z_{k+1},\xib)-f_r(z_{k+1},\xib)|&= \frac{1}{2}\left|\|y(z_{k+1},\xib)\|^2_{L^2(\D_0)}-\|y_r(z_{k+1},\xib)\|^2_{L^2(\D_0)}\right|\\
&\leq \frac{1}{2}\|y(z_{k+1},\xib)+y_r(z_{k+1},\xib)\|^2_{L^2(\D_0)}\|y(z_{k+1},\xib)-y_r(z_{k+1},\xib)\|^2_{L^2(\D_0)}\\
&\leq \frac{1}{2}\|y(z_{k+1},\xib)+y_r(z_{k+1},\xib)\|^2_{L^2(\D_0)}e^{\RB}_{z_{k+1}}(\xib).
\end{aligned}
\end{equation*}
In practice, as $y(z_{k+1},\xib)$ is not available, we use the estimator $e_{z_{k+1}}(\xib):=\|y_r(z_{k+1},\xib)\|^2_{L^2(\D)}e^{\RB}_{z_{k+1}}(\xib)$.

\subsection{Comparison of Alg. \ref{alg:adaptive_sampling} and Alg. \ref{alg:adaptive_sampling_rom}}

We now discuss numerical experiments aimed at comparing the standard adaptive-sampling algorithm (Alg \ref{alg:adaptive_sampling}) and its new variant leveraging importance sampling (Alg. \ref{alg:adaptive_sampling_rom}). All experiments are based on a $\PP^1$ finite element discretization of the state equation with 1023 degrees of freedom. The smoothing parameter for the $\Cvar$ is $\varepsilon=10^{-4}$. The optimization parameters are $\alpha=0.5$, $\theta=0.5$, $M_0=10$, $M_{trial}=10^3$, and all algorithms are stopped when a computational budget of $M_{\max}=10^6$ PDE solutions would be exceeded in the next iteration. For Alg. \ref{alg:adaptive_sampling_rom}, we do include the number of additional PDE solutions needed by the Greedy-enhanced construction of the ROM into the total count. The results shown are obtained starting from the initial guess $z_0=0$ and by averaging over 20 runs the randomness of the sample selection.

\begin{figure}
\centering
\includegraphics[width=0.33\linewidth]{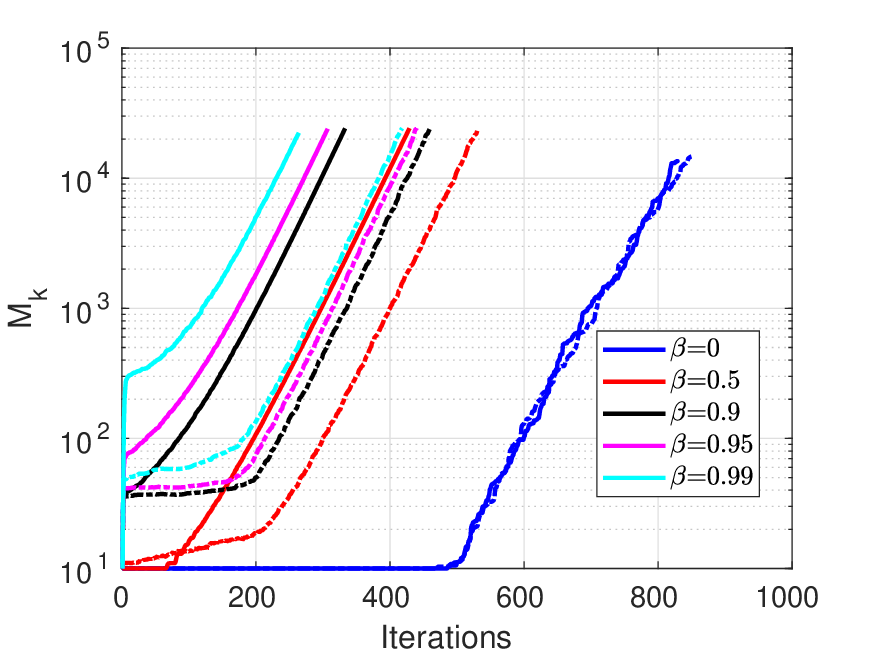}
\includegraphics[width=0.33\linewidth]{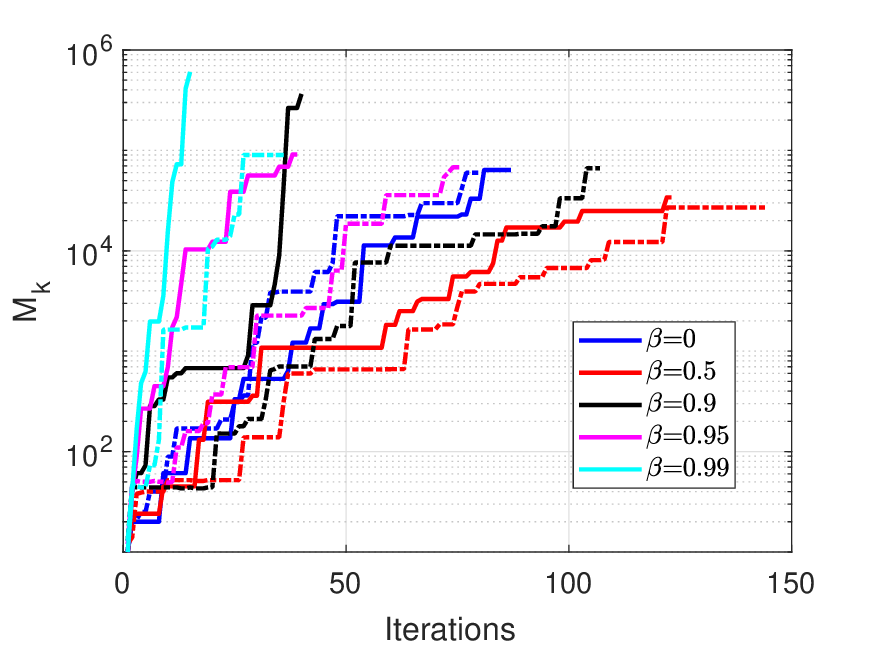}\\
\includegraphics[width=0.33\linewidth]{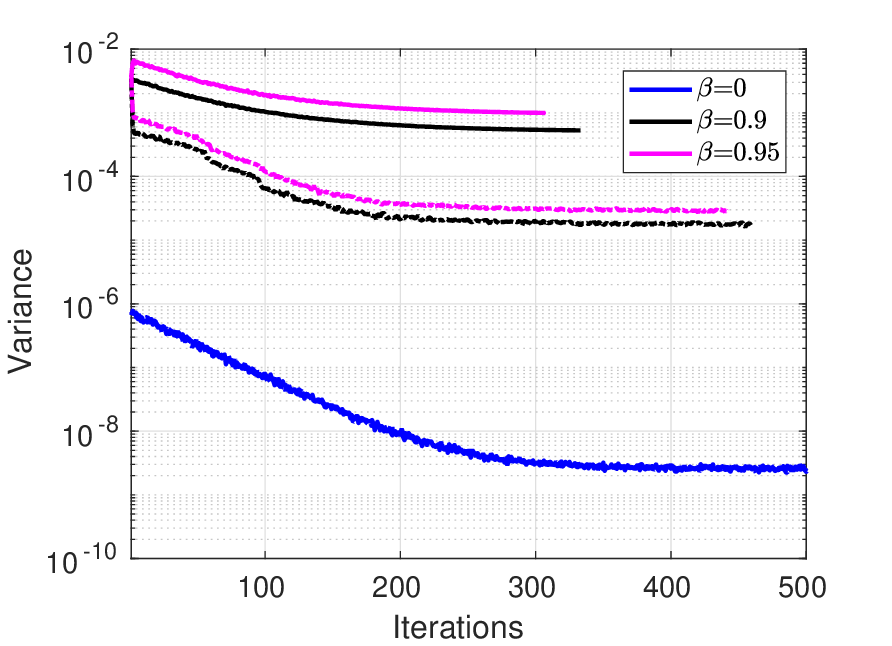}
\includegraphics[width=0.33\linewidth]{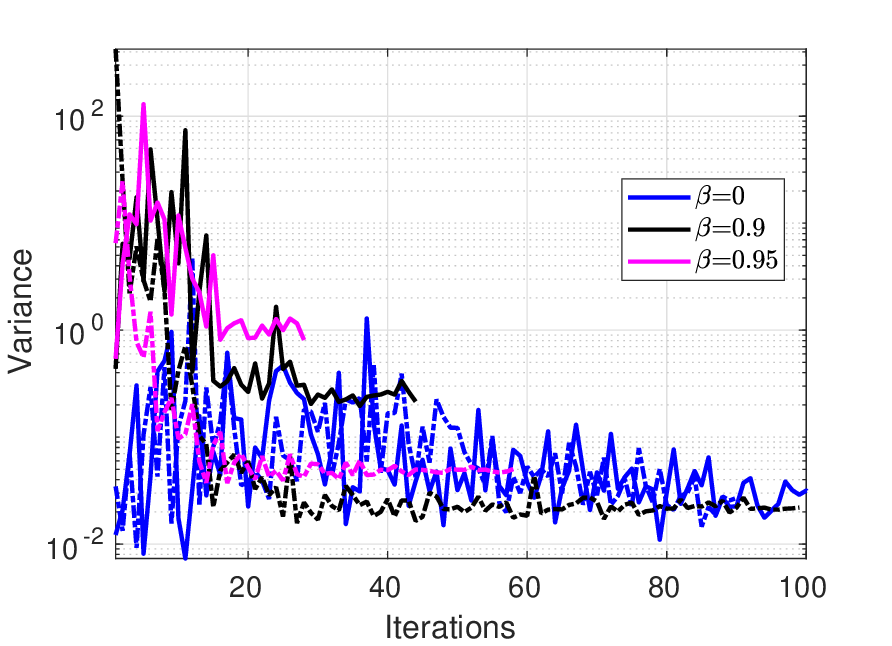}
\caption{Top row: growth of the sample size $M_k$ with iterations.
Bottom row: estimated empirical variance of the gradients along the iterations.
Left panels refer to $\kappa_1$, the right panels to $\kappa_2$. The continuous lines refer to Alg \ref{alg:adaptive_sampling} while the dashed-dotted ones to Alg. \ref{alg:adaptive_sampling_rom}.} \label{Fig:Mvec_beta}.  \end{figure}

The top plots of Figure \ref{Fig:Mvec_beta} show the growth of the sample sizes $M_k$ along the optimization for different values of $\beta$. The bottom plots show the empirical variance of the gradients that appears in the two norm tests \eqref{eq:norm_test_pra_cvar} and \eqref{eq:norm_test_cvar_pra_IS}.
As argued in Section \ref{sec:challenges}, larger values of $\beta$ lead to larger variances of the gradients, which in turn imply a faster growth of the sample sizes $M_k$.  Consequently, the limit of the computational budget is reached the earlier the larger is the value of $\beta$. By using importance sampling, we indeed achieve a variance reduction and hence Alg. \ref{alg:adaptive_sampling_rom} manage to perform additional iterations.

\begin{figure}
    \centering
    \includegraphics[width=0.3\linewidth]{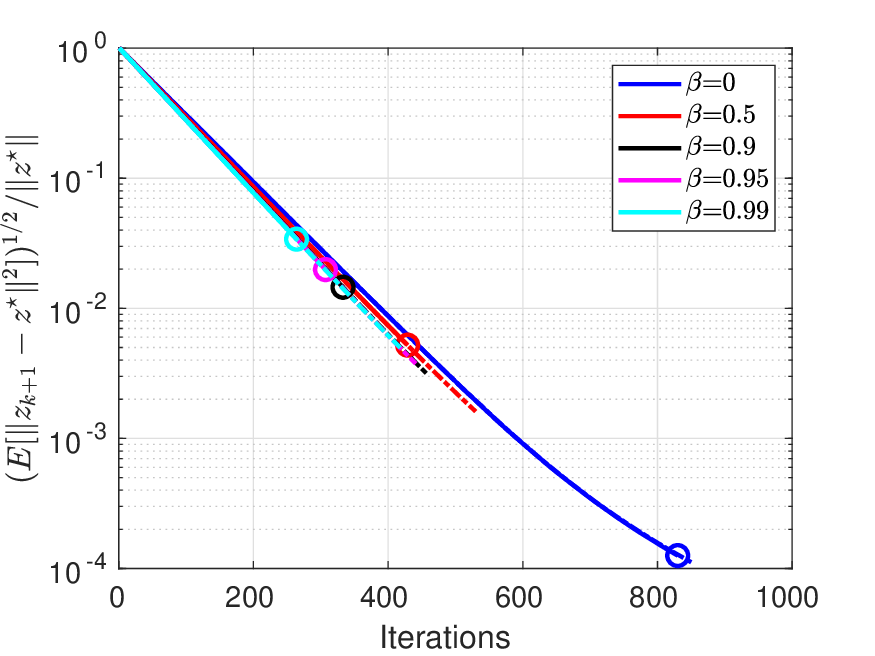}
    \includegraphics[width=0.3\linewidth]{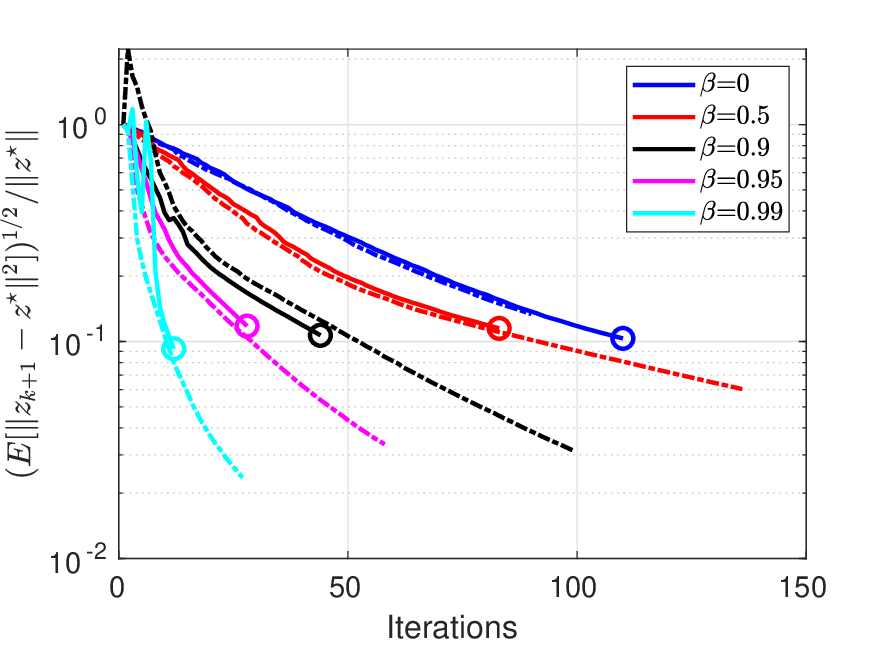}\\
   \includegraphics[width=0.3\linewidth]{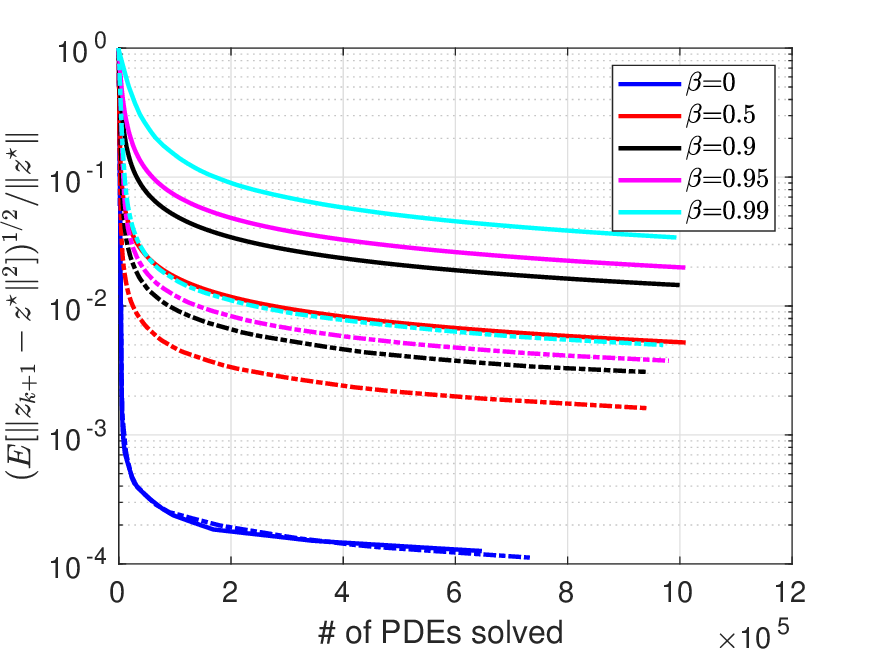}
    \includegraphics[width=0.3\linewidth]{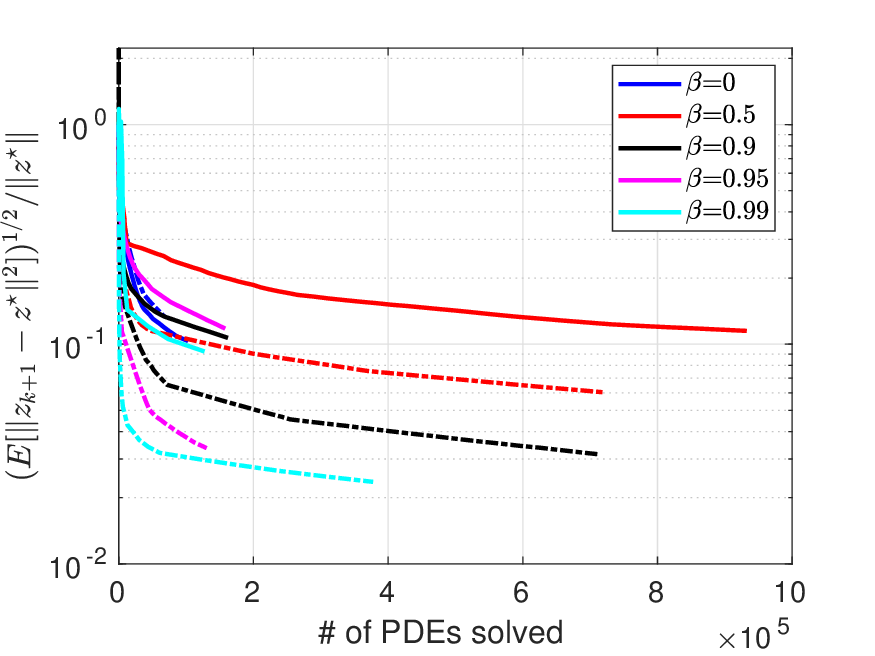}
  \caption{Decay of the relative error with respect to the number of iterations (top) and number of PDEs solved (bottom). The left column refers to $\kappa_1$, while the right one to $\kappa_2$. 
  The continuous line corresponds to Alg. \ref{alg:adaptive_sampling}, the dashed-dotted lines to Alg. \ref{alg:adaptive_sampling_rom}. The circle marker highlights the last iteration of Alg. \ref{alg:adaptive_sampling_rom}.} \label{Fig:Err_vs_it}
\end{figure}

Next, we consider the convergence behavior of both algorithms.
Figure \ref{Fig:Err_vs_it} shows the decay in expectation of the relative error with respect to the number of iterations and number of PDEs solved for different values of $\beta$. Since some lines are overlapping, we indicate with a circle marker the last iteration of Alg. \ref{alg:adaptive_sampling}.
The results are in agreement with the convergence analysis of Section \ref{sec:convergence} and with Theorem \ref{thm:convergence}. The asymptotic convergence rate is the same for both Alg. \ref{alg:adaptive_sampling} and Alg. \ref{alg:adaptive_sampling_rom}, but Alg. \ref{alg:adaptive_sampling_rom} can perform additional iterations, reaching significant smaller relative errors in all cases.

Finally, we compare the computational times. The results are obtained using a serial implementation, in which both the PDEs solves and the acceptance-rejection sampling are performed sequentially.
To fairly compare the algorithms, as they exhaust the computational budget at different iterations, Figure 4 reports the computational times required by both algorithms up to the iteration at which Alg. \ref{alg:adaptive_sampling} is stopped (circle marker in Figure \ref{Fig:Err_vs_it}). This is a reasonable metric, since Figure \ref{Fig:Err_vs_it} shows that both algorithms reach approximately the same relative error. Table \ref{tab:speed} reports the related speedup factors.
In all cases, Alg. \ref{alg:adaptive_sampling_rom} leads to an effective reduction of the computational times, sustaining the proposed methodology. 

\begin{figure}[h]
\centering
\includegraphics[width=0.3\linewidth]{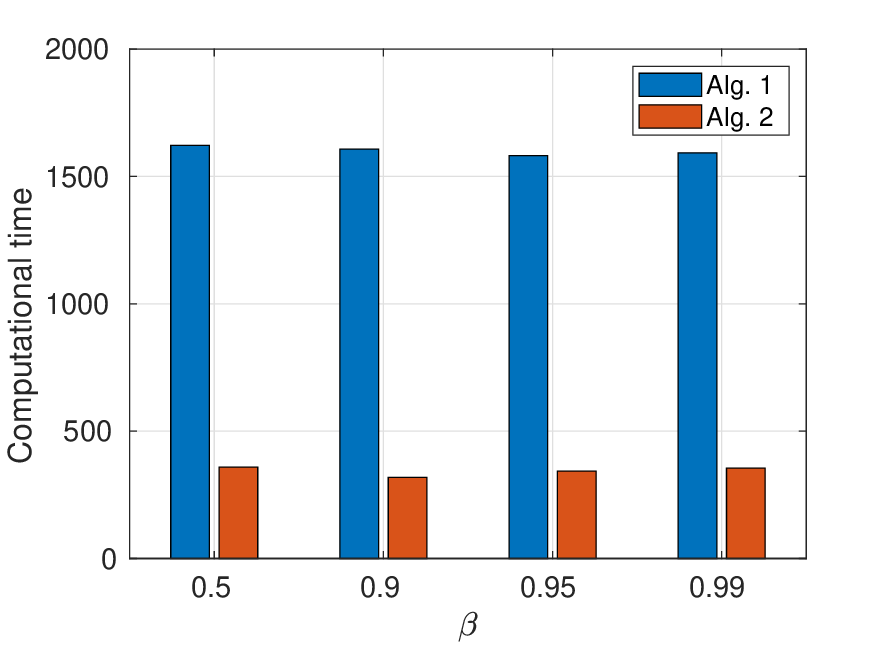}
\includegraphics[width=0.3\linewidth]{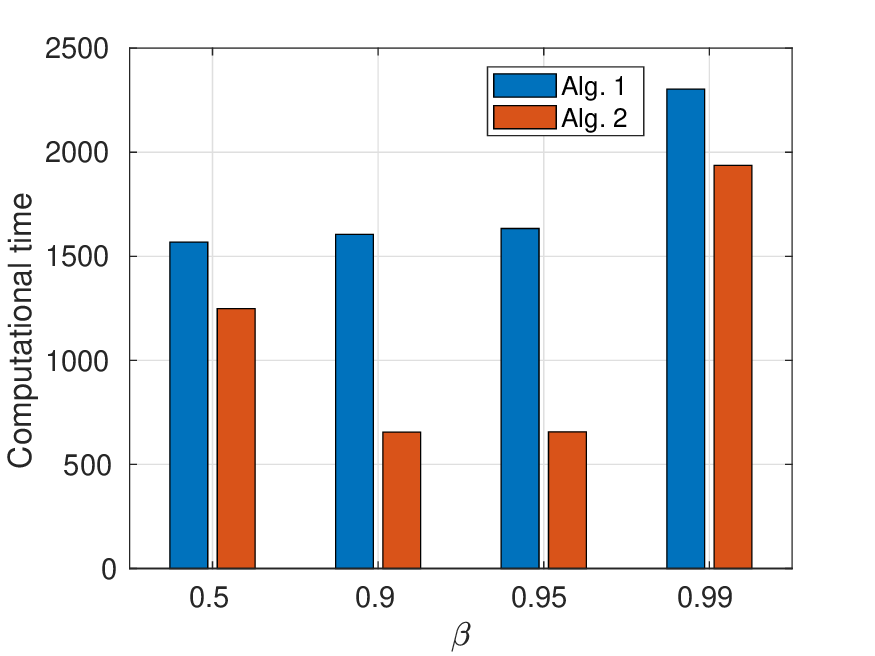}
\caption{Comparison of the computational times in seconds required by Alg. \ref{alg:adaptive_sampling} and Alg. \ref{alg:adaptive_sampling_rom} to reach approximately the same relative error. Left panel refers to $\kappa_1$, while the right one to $\kappa_2$.}
\end{figure}

\begin{table}[h]
\centering
\begin{tabular}{|c|c|c|c|c|}
\hline
$\beta$ & 0.5 & 0.9 & 0.95 & 0.99\\\hline
T\textsubscript{IS}/T  & 0.22 & 0.19 & 0.21 & 0.22 \\\hline
\end{tabular}
\begin{tabular}{|c|c|c|c|c|}
\hline
$\beta$ & 0.5 & 0.9 & 0.95 & 0.99\\\hline
T\textsubscript{IS}/T  & 0.79 & 0.4 & 0.4 & 0.84  \\\hline
\end{tabular}
\caption{Ratio between computational times.}\label{tab:speed}
\end{table}

\section{Conclusions}
This work has presented an improved adaptive sampling algorithm that adjusts both the sample size and the biasing distribution at each iteration of the optimization process. By restricting the sampling to a neighborhood of the risk region, the algorithm enables the use of smaller sample size per iteration, while preserving the linear asymptotic convergence rate.
As a consequence, the new methodology presents a better computational complexity measured in computational times, to achieve a desired relative error.

Future research could focus on further improving the on-the-fly construction of the reduced-order model (ROM). For instance, this could involve recycling the reduced basis subspace between iterations and augmenting it only when necessary. Another interesting direction would be the generation of the biasing distribution without relying on a ROM, by, for instance leveraging the information gathered from evaluating the full-order model at each iteration in combination with machine learning algorithms.

\section{Acknowledgements}
The authors are members of the INdAM-GNCS group. The first author acknowledges the support of the PRIN project 20227K44ME {\em Full and Reduced order modelling of coupled systems: focus on non-matching methods and automatic learning (FaReX)}. The second author has been partially supported by the INdAM-GNCS project {\em GNCS 2024 - CUP E53C23001670001}.

\section{Declaration of competing interest}
The authors declare that they have no known competing financial interests or personal relationships that could have appeared to
influence the work reported in this paper.

\section{Data availability}
No data was used for the research described in the article.

\section{CRediT authorship contribution statement}
\textbf{Sandra Pieraccini}:Writing -  review and editing, Methodology, Supervision. \textbf{Tommaso Vanzan}: Writing – review and editing, Writing- original draft, Software, Methodology, Formal analysis, Conceptualization.


\end{document}